# CONSISTENT COVARIATE SELECTION AND POST MODEL SELECTION INFERENCE IN SEMIPARAMETRIC REGRESSION

BY FLORENTINA BUNEA

*Florida State University*

This paper presents a model selection technique of estimation in semiparametric regression models of the type $Y_i = \beta'\underline{X}_i + f(T_i) + W_i$, $i = 1, \ldots, n$. The parametric and nonparametric components are estimated simultaneously by this procedure. Estimation is based on a collection of finite-dimensional models, using a penalized least squares criterion for selection. We show that by tailoring the penalty terms developed for nonparametric regression to semiparametric models, we can consistently estimate the subset of nonzero coefficients of the linear part. Moreover, the selected estimator of the linear component is asymptotically normal.

**1. Introduction.** The partially linear regression model was introduced by Engle, Granger, Rice and Weiss (1986) for the study of the relationship between weather and electricity sales and has received considerable attention over the last decade. Given $n$ i.i.d. $P$ observations $((\underline{X}_1, T_1), Y_1), \ldots, ((\underline{X}_n, T_n), Y_n)$, the model is

(1.1) $$Y_i = \beta'\underline{X}_i + f(T_i) + W_i = s(\underline{X}_i, T_i) + W_i,$$

where $W_1, \ldots, W_n$ are independent, identically distributed and zero mean error variables, assumed to be independent of $(\underline{X}, T) \in \mathbb{R}^q \times \mathbb{R}$. We assume that $f \in \mathcal{F}_\alpha$, where $\mathcal{F}_\alpha$ is a class of smooth functions with degree of smoothness $\alpha$. The appeal of the model lies in its flexibility. It can be used when a simple linear regression is adequate, apart from a covariate, usually a confounder, that is known to affect the response in a nonlinear fashion. More generally, (1.1) can be regarded as a particular case of a multiple index model and can serve as a first step in a dimension reduction process.

Two important aspects of optimality in estimating (1.1) are:

O1. Parsimonious selection of the $\underline{X} = \{X_1, \ldots, X_q\}$ covariates that avoids both underfitting and overfitting.
O2. Asymptotic normality of the *selected* estimator of $\beta$.

Let $I_0 \subseteq \{1, \ldots, q\}$ be the index set of the nonzero components of $\beta$. Let $\hat{I}$ be an estimator of $I_0$. We interpret $P(\hat{I} = I_0) \to 1$ as O1. The asymptotic normality of $\hat{\beta}_{\hat{I}}$, the estimator of $\beta$ corresponding to the selected index set $\hat{I}$, is our O2.

Model (1.1) was originally studied under the assumption that $I_0 = \{1, \ldots, q\}$ and $\alpha$ are known. We shall refer to this situation as Case 0. In this case, one can construct estimators of $f$ and $\beta$ using the knowledge of $\alpha$ and $I_0$, so there is no need for a model selection procedure. Wahba (1984), Green, Jennison and Seheult (1985) and Heckman (1986) suggested a least squares approach combined with spline smoothing for the nonlinear component; Chen (1988) proposed simultaneous least squares estimation of the parametric and nonparametric parts. The subsequent literature is vast, and we refer to Bickel, Klaassen, Ritov and Wellner (1993) and to the monograph by Härdle, Liang and Gao (2000) for an extensive bibliography. All suggested methods lead to asymptotically normal estimators of $\beta$, provided that the estimators of $f$ satisfy the minimum requirement

$$\|f - \hat{f}\|^2_{\mu_T} = o_P(n^{-1/2}), \tag{1.2}$$

for $\alpha > 1/2$, where $\|\cdot\|_{\mu_T}$ is the $L_2(\mu_T)$ norm and $\mu_T$ is the probability distribution of $T$; see, for instance, Lemma 11.2, page 202, in van der Geer (2000), or Chen (1988). Throughout this paper we shall use (1.2) as a prerequisite for O2 and assess it over functions in $\text{Lip}^*(\alpha, L_2(\mu_T))$, with $\alpha > 1/2$, defined in Section 2.1. We will compare the rate in (1.2) with the minimax rate in Corollary 3.1.

In this paper we use model selection based on penalized least squares minimization as an estimation procedure. We construct a sequence of finite-dimensional approximating spaces for $s$ and find the least squares estimator corresponding to each space. We compute the residual sum of squares corresponding to each estimator and then add a penalty term. Our final estimator is the one with the smallest penalized residual sum of squares. This yields estimators for $\beta$, $I_0$ and $f$ simultaneously.

The aim of this paper is to study the performance of such estimators when we relax the conditions of Case 0. We consider the following cases:

- Case 2: $I_0$ unknown and $\alpha$ known, $\alpha > 1/2$.
- Case 3: $I_0$ and $\alpha$ unknown, $\alpha > 1/2$.

We show in Sections 3.2 and 4 that O1 (the consistency of the selected index) and O2 (the asymptotic normality of the selected estimator) hold in both cases. We note that in Case 3 we need O1 and O2 to hold *uniformly* over the range of $\alpha$. Notice further that the smoothness classes are nested, in the sense that $\alpha_1 \geq \alpha_2$ implies $\text{Lip}^*(\alpha_1, L_2(\mu_T)) \subseteq \text{Lip}^*(\alpha_2, L_2(\mu_T))$. Then, to establish the uniformity result it is enough to show that O1 and O2 hold for the smallest allowable $\alpha > 1/2$. This is the approach we adopt in the sequel.

We point out that our strategy of showing O2, in either case, involves two steps. The first one is to show that $\hat{I}$ stabilizes asymptotically, that is, that O1 holds. The second one requires that the estimator of $\beta$ of fixed dimension $I_0$ be asymptotically normal. This could be done by quoting directly the result obtained in Case 0 by Chen (1988), whose method of estimation is the closest to ours. However, the result of Chen (1988) holds over $m$ times continuously differentiable functions $f$ (cf. his Condition 2, page 138), so it is not directly applicable to our context, which covers nondifferentiable functions. We therefore address this situation here, and in fact we prove that the asymptotic normality of $\hat{\beta}$ of fixed dimension holds in the more general case:

- Case 1: $I_0$ known and $\alpha$ unknown, $\alpha > 1/2$.

O2 was not studied in any of these three cases.

Partial results on O1 were obtained in Case 3. Härdle, Liang and Gao (2000), for a time series version of model (1.1), used a kernel method to estimate $f$ and cross validation to select $I_0$ and the bandwidth simultaneously. They show O1 in Theorem 6.3.1, page 137, but the construction of their estimator depends on the unknown $I_0$, as their Assumption 6.6.7, page 158, imposes lower and upper bound restrictions on the bandwidth that depend on $I_0$.

Chen and Chen (1991) used $B$-splines to approximate the nonlinear component and an Akaike type technique for simultaneous estimation of $I_0$ and $f$. They discuss O1 in their Proposition 2, page 334, under a condition on a random criterion that depends on the unknown $I_0$ and $f$. Remark 3 verifies this condition only for $|I_0| = 1$ and under their Condition 4, page 326, which entails the existence of a lower bound on $f$; the further study of O1 is left as an open problem. We do not use any of these conditions here.

The remainder of the paper is organized as follows. Section 2.1 contains the assumptions under which our results hold. In Section 2.2 we give the construction of our estimators. In Section 2.3 we derive upper-bound oracle inequalities for the risk of our estimators, and obtain rates of convergence as a consequence. Sections 3 and 4 are central to our paper. In Section 3.1 we discuss penalty choices and their impact on the rates of convergence derived in Section 2.3. We prove O1 in Section 3.2. In Section 4 we show that the estimators of $\beta$ are asymptotically normal for the three cases under consideration, therefore establishing O2. Section 5 provides conclusions. The proofs of intermediate results are given in the Appendix.

**2. A penalized least squares estimator.** In this section we devise estimators for $\beta$ and $f$ and establish their consistency.

2.1. *Preliminaries.* We begin by giving a list of assumptions under which the results of this paper hold.

Let $\mu$ be the joint probability distribution of $\underline{X}$ and $T$. Let $\mu_X$ and $\mu_T$ denote the probability distributions of $\underline{X}$ and $T$, respectively. Let $\text{Lip}^*(\alpha, L_2(\mu_T))$ denote a generalized Lipschitz space, for some smoothness parameter $\alpha > 0$, and let $|\cdot|_{\alpha,2}$ be the seminorm in this space [see, e.g., DeVore and Lorentz (1993), page 51, for definition and properties]. For some positive constant $A > 0$ define $\mathcal{D}_{\alpha,2}(A) = \{g \in \text{Lip}^*(\alpha, L_2(\mu_T)), |g|_{\alpha,2} \leq A\}$.

ASSUMPTION 2.1. The support of $\mu$ is $[0,1] \times \mathcal{K}$, for some compact set $\mathcal{K} \subset \mathbb{R}^q$, and on its support $\mu$ admits a density with respect to the Lebesgue measure $\lambda$ that is bounded from below by $h_0 > 0$ and from above by $h_1 < \infty$. In addition, $\mu_T([0,1]) = 1$ and $\mu_X(\mathcal{K}) = 1$.

ASSUMPTION 2.2. There exists $d > 0$ such that $\tau_p = E(|W_1|^p) < \infty$ for $p > 4 + d$.

ASSUMPTION 2.3. For any $l \in \mathbb{R}^q \setminus \{0\}$, $\text{Var}(l'\underline{X}|T = t) > 0$ for all $t \in [0,1]$.

ASSUMPTION 2.4. $s \in L_2(\mathcal{K} \times [0,1], \lambda)$.

ASSUMPTION 2.5. $\theta_j = E(X_j|T = t) \in \mathcal{D}_{\gamma,2}(A)$ for some $\gamma > b/4$ and some fixed constant $b \geq 3$, for all $j \in \{1, \ldots, q\}$.

Note that Assumption 2.3 ensures that model (1.1) is identifiable; see, for example, Lemma 11.2 of van der Geer (2000) or Lemma 3 of Chen (1988).

Assumption 2.5 is a sufficient condition on the smoothness of $\theta_j$, $j \in \{1, \ldots, q\}$, which ensures that $\beta$ can be estimated at the optimal $n^{-1/2}$ convergence rate. See also Chen (1988), Heckman (1986), Speckman (1988) and van der Geer (2000) for other types of smoothness conditions, typically requiring the existence of a prespecified number of derivatives for $\theta_j$.

2.2. *The sieves and the estimators.* We construct now a sequence of approximating spaces for $s$ in (1.1). This mimics the construction of approximating spaces for generalized additive models, as in Barron, Birgé and Massart (1999) or Baraud (2000). Let $I = \{i_1, \ldots, i_l\} \in \mathcal{I}$, where $\mathcal{I} = \mathcal{P}(\{1, \ldots, q\})$, and $\mathcal{P}(F)$ denotes all subsets of a set $F$. Denote by $[\cdot]$ the integer part and by $\log_2$ the logarithm in base 2. Let $b$ be a fixed constant, $b \geq 3$. Define

$$
\begin{aligned}
A_n &= [\log_2(n/\log n)^{1/b}], \\
J_n &= [\log_2(n/\log n)^{1/2}], \\
B_n &= 2^{A_n}, \qquad N_n = 2^{J_n}.
\end{aligned}
\tag{2.1}
$$

For $k_n \in \{A_n, A_n+1, \ldots, J_n\}$ let $K_n = 2^{k_n}$. For each $K_n \in \{B_n, \ldots, N_n\} = \mathcal{K}_n$ let $S_{K_n}$ be the linear space of piecewise polynomials of degree at most $r-1$, based on a regular dyadic partition of size $1/K_n$. Thus, $S_{K_n}$ is the space of functions $v$ on $[0,1]$ of the form $v(t) = \sum_{j=1}^{K_n} P_j(t)\mathbb{1}(\{\frac{j-1}{K_n} \leq t < \frac{j}{K_n}\})$, where $\mathbb{1}(V)$ denotes the indicator of a set $V$. Note that $\dim(S_{K_n}) = rK_n$. Denoting the restriction of $\lambda$ to $[0,1]$ by $\lambda_T$, let $\{\phi_j\}_{j=1}^{rK_n}$ be an orthonormal basis in $L_2(\lambda_T)$ for $S_{K_n}$. For $(I, K_n) \in \mathcal{I} \times \mathcal{K}_n$ define

$$S_{I,K_n} = \langle x_{i_1}, \ldots, x_{i_l}, \phi_1(t), \ldots, \phi_{rK_n}(t)\rangle,$$

where $\langle \cdot \rangle$ denotes the linear span. Note that $|\mathcal{I} \times \mathcal{K}_n| = 2^q \times (J_n - A_n + 1)$, where $|\cdot|$ denotes the cardinality of a set. Notice that $\dim(S_{I,K_n}) = |I| + rK_n$.

We recall that the approximating space $S_{K_n}$ is known to have good approximation properties for a range of smoothness classes to which $f$ and $\theta_j$, $j = 1, \ldots, q$, may belong; see, for example, DeVore and Lorentz (1993). However, other choices are possible: spaces generated by piecewise polynomials based on an irregular partition of $[0,1]$, wavelets or trigonometric polynomials; see Birgé and Massart (1998) for a detailed discussion.

Let $\text{pen}(I, K_n)$ be a penalty term associated with $S_{I,K_n}$. We defer a detailed discussion on the penalty term to Section 3.1. For $(I, K_n) \in \mathcal{I} \times \mathcal{K}_n$ and $u \in S_{I,K_n}$ let $\gamma_n(u) = n^{-1}\sum_{i=1}^{n}[Y_i - u(\underline{X}_i, T_i)]^2$.

DEFINITION 2.1. A penalized least squares estimator relative to the collection $\{S_{I,K_n}\}_{(I,K_n) \in \mathcal{I} \times \mathcal{K}_n}$ is any $\tilde{s} \in S_{\hat{I}, \hat{K}_n}$ such that

$$\gamma_n(\tilde{s}) + \text{pen}(\hat{I}, \hat{K}_n) = \inf_{(I, K_n) \in \mathcal{I} \times \mathcal{K}_n}\left(\inf_{u \in S_{I,K_n}} \gamma_n(u) + \text{pen}(I, K_n)\right). \tag{2.2}$$

Let $\mathbf{Y} = (Y_1, \ldots, Y_n)'$. Denote by $\mathbf{X}$ the $n \times q$ matrix with columns $(X_{1,j}, \ldots, X_{n,j})'$, $1 \leq j \leq q$, and by $\mathbf{X}_I$ the $n \times |I|$ matrix obtained from $\mathbf{X}$ by retaining the columns corresponding to the index set $I \subseteq \{1, \ldots, q\}$. Let $\beta_I$ be a vector in $\mathbb{R}^{|I|}$ and let $\delta_{K_n}$ be a vector in $\mathbb{R}^{rK_n}$. Let $\mathbf{Z}_{K_n}$ be the $n \times rK_n$ matrix whose $i$th row is $\phi_1(T_i), \ldots, \phi_{rK_n}(T_i)$. Then, in matrix notation, our estimator achieves the infimum below:

$$\inf_{(I,K_n)} \inf_{\beta_I, \delta_{K_n}} \{(\mathbf{Y} - \mathbf{X}_I\beta_I - \mathbf{Z}_{K_n}\delta_{K_n})'(\mathbf{Y} - \mathbf{X}_I\beta_I - \mathbf{Z}_{K_n}\delta_{K_n}) + \text{pen}(I, K_n)\}. \tag{2.3}$$

Let $\hat{K}_n$ and $\hat{I}$ be the indices for which (2.3) is attained. Then, if the minimization problem has a unique solution, following Seber [(1977), Theorem 3.7], the least squares estimators of $\beta_I$ and $\delta_{K_n}$ are, respectively,

$$\tilde{\beta}_{\hat{I}} = (\mathbf{X}'_{\hat{I}}(\mathbf{Id} - \mathbf{P}_{\hat{K}_n})\mathbf{X}_{\hat{I}})^{-1}\mathbf{X}'_{\hat{I}}(\mathbf{Id} - \mathbf{P}_{\hat{K}_n})\mathbf{Y} \tag{2.4}$$

and

$$\tilde{\delta}_{\hat{K}_n} = (\mathbf{Z}'_{\hat{K}_n}\mathbf{Z}_{\hat{K}_n})^{-1}\mathbf{Z}'_{\hat{K}_n}(\mathbf{Y} - \mathbf{X}_{\hat{I}}\tilde{\beta}_{\hat{I}}), \tag{2.5}$$

where $\mathbf{P}_{\hat{K}_n}$ is the projection matrix on the space $L_{\hat{K}_n}$ generated by the columns of $\mathbf{Z}_{\hat{K}_n}$. Thus, $L_{\hat{K}_n} = \{(g(T_1), \ldots, g(T_n))' | g \in S_{\hat{K}_n}\}$ and $\mathbf{P}_{\hat{K}_n} = \mathbf{Z}_{\hat{K}_n}(\mathbf{Z}'_{\hat{K}_n} \times \mathbf{Z}_{\hat{K}_n})^{-1}\mathbf{Z}'_{\hat{K}_n}$. $\mathbf{Id}$ denotes the $n \times n$ identity matrix.

For any measure $\nu$ we denote by $\|\cdot\|_\nu$ the $L_2(\nu)$ norm. Let

$$\Gamma_n = \{\|\tilde{s}\|_\lambda \leq 2\exp(\log^2 n)\}.$$

We note that, by Theorem 1.1 in Baraud (2002) $P(\Gamma_n^c) \to 0$. For technical reasons we consider as our final estimator of $s$,

$$\hat{s}(x,t) = \tilde{s}(x,t)\mathbb{1}_{\Gamma_n}.$$

We denote by $\tilde{f}$ the estimator of $f$ corresponding to $\tilde{\delta}_{\hat{K}_n}$. Hence, the estimators of the nonlinear and linear part are, respectively,

$$\hat{f} = \tilde{f}\mathbb{1}_{\Gamma_n} \quad \text{and} \quad \hat{\beta} = \tilde{\beta}\mathbb{1}_{\Gamma_n}.$$

We mention here the approximating spaces used in the three cases under consideration. We elaborate on this in Section 3.1. For Case 1 we use $\{S_{I_0,K_n}\}_{K_n \in \mathcal{K}_n}$. In Case 2 we use $\{S_{I,K_{n,\alpha}}\}_{I \in \mathcal{I}}$ with $K_{n,\alpha} \asymp n^{1/2\alpha+1}$. Here and in the sequel the notation $a \asymp b$ means that $a$ is an integer power of 2 that differs from $b$ by at most a factor of 2. In Case 3 we use $\{S_{I,K_{n,a}}\}_{I \in \mathcal{I}}$ for $K_{n,a} \asymp n^{1/2a+2}$, with $a > 0$ arbitrarily close to zero. Notice that $K_{n,\alpha} \in \mathcal{K}_n$ if $1/2 < \alpha < (b-1)/2$. We shall need later the approximation theory result given by (2.7), which holds for $\alpha \in (0,r)$. This motivates the choice of $b \geq 3$ and the definition of $r = [\frac{b-1}{2}]$. Also, note that $K_{n,a} \asymp n^{1/2a+2} \in \mathcal{K}_n$ for any $0 < a < 1/2$.

We show in Appendix A.1 that, under Assumptions 2.1–2.5, the estimators are unique, except for a set of probability tending to zero. On this set we define our estimators to be identically zero.

2.3. *The consistency of the penalized least squares estimators.* In this section we give finite sample upper bounds on the risk of the estimator $\hat{s}$ of $s$. As an immediate consequence we then obtain rates of convergence for the estimators of $\beta$ and $f$, respectively. Oracle type inequalities for the risk of estimators obtained via model selection in nonparametric regression have been studied extensively over the last decade; see, for example,

Barron, Birgé and Massart (1999), Baraud (2000, 2002), Wegkamp (2003) and the references therein. The results carry over to semiparametric regression and in this paper we adopt the approach of the second author.

Let $f_{K_n}$ be the $L_2(\mu_T)$ projection of $f$ onto $S_{K_n}$. Let $C_1, C_2 > 0$ denote dominating constants independent of $n$, given by Theorem 2.1 in Baraud (2002).

THEOREM 2.1. *Under Assumptions 2.1–2.5, for* $\mathrm{pen}(I, K_n) \geq C_1(|I| + rK_n)/n$,

$$(2.6) \qquad E\|s - \hat{s}\|_\mu^2 \leq C_2 \inf_{K_n \in \mathcal{K}_n} \{\|f - f_{K_n}\|_{\mu_T}^2 + \mathrm{pen}(I_0, K_n) + \vartheta_n\},$$

*with* $\vartheta_n = 1/n + (\|s\|_\lambda^2 + 1)\exp(-2\log^2 n)$.

We note that $C_1$ depends on $\tau_2$ and that $C_2$ depends on $h_0, h_1, d, p, \tau_p$.

This theorem allows us to obtain rates of convergence for $\hat{s}$ by computing the infimum above, provided that $\|s\|_\lambda^2$ is bounded, in which case $\vartheta_n = O(n^{-1})$. This is guaranteed under Assumption 2.1 if $f \in \mathcal{D}_{\alpha,2}(L)$ and $\beta \in \mathcal{K}_1$, for some compact set $\mathcal{K}_1 \in \mathbb{R}^q$. Furthermore, the next corollary shows that the rate of convergence of $\hat{s}$ is inherited by the estimators of $\beta$ and $f$. Note that if $f \in \mathcal{D}_{\alpha,2}(L)$ for some $L > 0$, then by Theorem 2.4, page 358, in DeVore and Lorentz (1993) and Assumption 2.1, for any $\alpha \in (0, r)$ there exists $C(\alpha, L) > 0$ such that

$$(2.7) \qquad \|f - f_{K_n}\|_{\mu_T}^2 \leq h_1 C(\alpha, L) K_n^{-2\alpha}.$$

Define

$$(2.8) \qquad r_n = \inf_{\mathcal{K}_n}\{h_1 C(\alpha, L) K_n^{-2\alpha} + \mathrm{pen}(I_0, K_n) + \vartheta_n\}.$$

Let $|\cdot|_2$ denote the Euclidean norm. For a function $g$ of generic argument $Z$ we denote its empirical norm by $\|g\|_n^2 = n^{-1}\sum_{i=1}^n g^2(Z_i)$. In addition, by abuse of notation, we regard here $\hat{\beta}_{\hat{I}}$ as a vector in $\mathbb{R}^q$ by adding 0's to the necessary positions.

COROLLARY 2.1. *Under Assumptions 2.1–2.5, for* $\mathrm{pen}(I, K_n) \geq C_1(|I| + rK_n)/n$, *if* $f \in \mathcal{D}_{\alpha,2}(L)$, $0 < \alpha < r$, *we have:*

1. $\|\hat{f} - f\|_{\mu_T}^2 = O_P(r_n)$, *uniformly over* $f \in \mathcal{D}_{\alpha,2}(L)$ *and* $\beta \in \mathcal{K}_1$.
2. $\|\hat{f} - f\|_n^2 = O_P(r_n)$, *uniformly over* $f \in \mathcal{D}_{\alpha,2}(L)$ *and* $\beta \in \mathcal{K}_1$.
3. $|\hat{\beta}_{\hat{I}} - \beta|_2^2 = O_P(r_n)$, *uniformly in* $\beta \in \mathcal{K}_1$.

We present the proof of these two results in Appendix A.2. We discuss in the next section our penalty choices and the corresponding values of $r_n$.

Although $\hat{f}$ may achieve the minimax optimal rate of convergence, Corollary 2.1 gives suboptimal rates for the estimator of $\beta$. We show in Section 4 that we can achieve the expected $n^{-1/2}$ rate of convergence for $\hat{\beta}_{\hat{I}}$, provided that $P(\hat{I} = I_0) \to 1$. In the next section we discuss penalty choices for which this holds.

**3. Penalty choices and the consistency of the selected index.** There has been a vast literature on the estimation of $I_0$ in the fully parametric context and we only mention here the seminal works of Mallows (1973), Akaike (1974), Schwarz (1978) and Shibata (1981). Typically, model selection procedures based on a penalized criterion use penalty terms that are proportional either to $|I|/n$, where $|I|$ is the dimension of a fitted model, or to $|I|\log n/n$. These give rise to Akaike-type (AIC) and Schwarz-type (BIC) methods, respectively. In AIC the selected model is expected to include about one superfluous parameter [Woodroofe (1982)]. BIC chooses the correct model with probability converging to 1 [Haughton (1988)]. See also Guyon and Yao (1999) for a recent survey.

In semiparametric models we cannot obtain the consistency of $\hat{I}$ by a simple extension of these methods, because a penalty term proportional to $(|I| + K_n)\log n/n$ no longer suffices. In the parametric case a penalty term essentially balances out the residual variance of competing models of different dimensions, and the bias term disappears for the models that include the true one. However, in the semiparametric case the penalty term is also required to balance out the bias introduced by approximating $f$ within a finite-dimensional space. Since in general $f$ does not belong to any of the approximating spaces, this bias is not zero, and a penalty that is proportional to the dimension of a fitted model is too small to achieve the correct balance.

3.1. *Penalty choices and rates of convergence.* In this section we give sufficient conditions on the penalty term for which the optimality criteria O1 and O2 hold.

We first discuss O1, which is $P(\hat{I} \neq I_0) \to 0$. Note that

(3.1) $$P(\hat{I} \neq I_0) = P(I_0 \not\subset \hat{I}) + P(I_0 \subsetneq \hat{I}).$$

We showed in Corollary 2.1 that the estimators of $\beta$ are consistent, for any penalty term that satisfies $\text{pen}(I, K_n) \geq (|I| + rK_n)\log n/n$. We will show in Theorem 3.1 that the first term in (3.1) converges to zero, for any penalty term that satisfies this restriction. However, we cannot use the consistency of $\hat{\beta}_{\hat{I}}$ to show that the second term converges to zero, as we can overestimate the model but still consistently estimate 0's. The study of the convergence to zero of the second term in (3.1) leads to the second set of restrictions on our penalty term, namely, $\text{pen}(I, K_n) - \text{pen}(I_0, K_n) > h_1 C(\alpha) L K_n^{-2\alpha}$. This condition means that the penalty term needs to be greater than the bias induced by approximating $f$. Intuitively, if the bias due to the nonparametric component is present, it acts as a confounder, and the true parametric dimension cannot be recovered. Formally, as in (3.21), this condition ensures that $I_0$ can be found asymptotically.

We show in Theorem 3.1 below that O1 holds if the two conditions below are satisfied simultaneously:

(3.2) $\quad$ (i) $\text{pen}(I, K_n) - \text{pen}(I_0, K_n) > h_1 C(\alpha) L K_n^{-2\alpha}$,
$\quad\quad\,\,\,$ (ii) $\text{pen}(I, K_n) \geq (|I| + r K_n) \log n / n$.

We discuss now sufficient conditions on the penalty term under which O2 holds. Recall that (1.2) is a prerequisite for O2. With the notation of Section 2.3, (1.2) can be rewritten as

(3.3) $\quad \sqrt{n} r_n \to 0 \qquad \text{for } r_n = \inf_{\mathcal{K}_n}\{h_1 C(\alpha, L) K_n^{-2\alpha} + \text{pen}(I_0, K_n) + \vartheta_n\}$.

We show in Theorems 3.1 and 4.2 that O1 and O2 are compatible if (3.2) and (3.3) hold simultaneously. Note now the apparent contradiction between (3.2)(i) and (3.3): the first one essentially requires that the penalty term dominate the bias for *all* $K_n$, whereas the oracle inequality of Theorem 2.1 tells us that the best rate of convergence of $\hat{f}$ is achieved for *a particular* $K_n$, namely, the one realizing the best bias-variance trade-off. Notice further that by simply taking $\text{pen}(I, K_n) = (|I| + r K_n) \log n / n$, we would have (3.2)(i) satisfied uniformly over $\alpha \in (1/2, r)$ only if, up to multiplicative constants, $K_n > n / \log n$ for *all* $K_n \in \mathcal{K}_n$, in which case the estimators would no longer be defined.

The first part of the solution is to construct a penalty term in which the dimensions $|I|$ and $rK_n$ are multiplied rather than added. Thus, we first consider

(3.4) $\qquad\qquad \text{pen}(I, K_n) = 2(|I| + 1) r K_n \log n / n$.

This penalty satisfies (3.2)(ii) and we show in Corollary 3.1 that it also leads to an $r_n$ that satisfies (3.3).

However, if we use (3.4) for either Case 2 or Case 3, then (3.2)(i) holds only if, up to multiplicative constants, $K_n > (n/\log n)^{1/2\alpha+1}$ for all $K_n \in \mathcal{K}_n$.

If $\alpha > 1/2$ is known, then for $K_{n,\alpha} \asymp n^{1/2\alpha+1}$ and $n$ large enough,

(3.5) $\qquad\qquad \text{pen}(I, K_{n,\alpha}) = 2(|I| + 1) r K_{n,\alpha} \log n / n$

satisfies (3.2) by construction, and (3.3) holds by Corollary 3.1, *for this* $\alpha$. We use the penalty term (3.5) in Case 2, for the approximating spaces $\{S_{I, K_{n,\alpha}}\}_{I \in \mathcal{I}}$.

If $\alpha$ is not known, as in Case 3, we can no longer define a penalty term that depends on $\alpha$. In this case we need to construct a penalty that satisfies the contradictory requirements (3.2)(i) and (3.3) for all $K_n$ and uniformly over $\alpha > 1/2$. Since they cannot hold simultaneously *for all $K_n$*, the strategy we suggest is to find the *best $K_n$* for which they are met, uniformly over $\alpha > 1/2$. For this, first recall that the smoothness spaces $\text{Lip}^*(\alpha, L_2(\mu_T))$ are nested: the smaller the $\alpha$, the larger the space and, also, the smaller the $\alpha$, the larger the bias term $1/K_n^\alpha$. Thus, the largest bias that needs to be dominated by the penalty term will correspond to $\alpha = 1/2 + a$ for $a > 0$ and arbitrarily close to zero, since this is the worst allowable case of $\alpha$. This reduces the penalty choice to the one of Case 2, for this particular choice of $\alpha$; namely, we choose $K_{n,a} \asymp n^{1/2a+2}$ and

$$(3.6) \qquad \text{pen}(I, K_{n,a}) = 2(|I| + 1)rK_{n,a} \log n/n.$$

The corresponding approximating spaces are in this case $\{S_{I,K_{n,a}}\}_{I \in \mathcal{I}}$. Then, *uniformly over $\alpha$*, (3.2) holds by construction, and (3.3) holds by Corollary 3.1.

We remark now that in Case 1, in which $I_0$ is considered known, (3.2)(ii) is no longer required. Only (3.2)(i), which also guarantees the consistency of the estimator of $\beta$, and (3.3) are needed. We shall therefore use the penalty term (3.4), in connection with the approximating spaces $\{S_{I_0,K_n}\}_{K_n \in \mathcal{K}_n}$, for this case.

We conclude this section with Corollary 3.1, which is an immediate consequence of Theorem 2.1. This result summarizes the rates of convergence corresponding to each penalty term, and shows that (1.2) holds in each case. We give the proof in Appendix A.2.

COROLLARY 3.1. *Under Assumptions* 2.1–2.5, *if $f \in \mathcal{D}_{\alpha,2}(L)$, we have:*

1. $\alpha$ *unknown.*

    (a) *If* $\text{pen}(I_0, K_n) = 2(|I_0| + 1)rK_n \log n/n$, *then*

    $$r_n = O((\log n/n)^{2\alpha/2\alpha+1}),$$

    *for any $\alpha \in (1/2, r)$ and known $I_0$.*
    (b) *If* $\text{pen}(I, K_n) = 2(|I| + 1)rK_{n,a} \log n/n$, *then*

    $$r_n = O(\log n/n^{(2a+1)/(2a+2)}),$$

    *for any $\alpha \in [1/2 + a, r)$, with $a > 0$ arbitrarily close to zero.*

2. $\alpha$ *known.*
    *If* $\text{pen}(I, K_{n,\alpha}) = 2(|I| + 1)rK_{n,\alpha} \log n/n$ *and $\alpha \in (1/2, r)$, then*

    $$r_n = O(\log n/n^{2\alpha/2\alpha+1}).$$

REMARK 3.1. Notice that in 1(a) of the above corollary $r_n$ is the minimax adaptive rate of convergence, up to a $\log n$ factor. Since $2\alpha/2\alpha + 1 > 1/2$ for any $\alpha > 1/2$, (1.2) holds uniformly over $\alpha \in (1/2, r)$.

For part 2, $r_n$ is the minimax rate, up to a $\log n$ factor, and (1.2) holds for the particular fixed $\alpha$.

For 1(b), although the rate is suboptimal relative to minimax, we have $(2a+1)/(2a+2) > 1/2$ for any $a > 0$, and so again (1.2) holds uniformly over $\alpha \in [1/2 + a, r)$.

REMARK 3.2. We note that the rate of convergence in 1(a) is achieved for $K_n^* \asymp (n/\log n)^{1/2\alpha+1}$. This dimension belongs to our $\mathcal{K}_n$ for all $1/2 < \alpha < r = [(b-1)/2]$.

### 3.2. The consistency of the selected index.

THEOREM 3.1. If $f \in \mathcal{D}_{\alpha,2}(L)$ and Assumptions 2.1–2.5 hold, then:

(a) $P(\hat{I} \neq I_0) \to 0$ as $n \to \infty$, for the penalty term given by (3.5) and for some given $\alpha \in (1/2, r)$.

(b) $P(\hat{I} \neq I_0) \to 0$ as $n \to \infty$, for the penalty term (3.6), uniformly over $\alpha \in [1/2 + a, r)$ for some $a > 0$, arbitrarily small.

PROOF. Notice that

$$P(\hat{I} \neq I_0) = P(I_0 \subsetneq \hat{I}) + P(I_0 \not\subset \hat{I}). \tag{3.7}$$

We show that each term in (3.7) converges to zero.

1. $P(I_0 \subsetneq \hat{I}) \to 0$ as $n \to \infty$.

(a) Notice that if $I_0 = \{1, \ldots, q\}$, $P(I_0 \subsetneq \hat{I}) = 0$, so it is enough to consider $I_0 \subsetneq \{1, \ldots, q\}$. Then

$$\lim_{n \to \infty} P(I_0 \subsetneq \hat{I}) = \lim_{n \to \infty} \sum_{I \supset I_0} P(\hat{I} = I). \tag{3.8}$$

For $I \supset I_0$ and re-denoting $K_n = K_{n,\alpha} \asymp n^{1/2\alpha+1}$, define

$$f_n(I, K_n) = \inf_{v \in S_{I,K_n}} \{\gamma_n(v) + \operatorname{pen}(I, K_n)\}. \tag{3.9}$$

Recall that $f_{K_n}$ is the $L_2(\mu_T)$ projection of $f$ onto $S_{K_n}$. Define $s_n(\underline{x}, t) = \beta_{I_0}\underline{x} + f_{K_n}(t)$, with $\beta_{I_0}$ regarded as a vector in $\mathbb{R}^q$, by adding zero to the necessary positions. Notice that $s_n \in S_{I_0, K_n}$ and so $s_n \in S_{I, K_n}$ for all $I \supset I_0$. Also, by (3.9) note that

$$f_n(I_0, K_n) \leq \gamma_n(s_n) + \operatorname{pen}(I_0, K_n).$$

Then, by the definition of the estimator, we have

$$P(\hat{I}=I) = P(f_n(I,K_n) - f_n(I',K'_n) \leq 0, \text{ for all } I' \neq I)$$

(3.10)
$$\leq P(f_n(I,K_n) - f_n(I_0,K_n) \leq 0)$$

$$\leq P\bigg(\sup_{v \in S_{I,K_n}} (\gamma_n(s_n) - \gamma_n(v)) \geq \text{pen}(I,K_n) - \text{pen}(I_0,K_n)\bigg).$$

For any function $g$ of generic argument $Z$ we denote $\mathbf{g} = (g(Z_1),\ldots,g(Z_n))'$. For any $\mathbf{U} = (U_1,\ldots,U_n)'$ we let $\|\mathbf{U}\|_n^2 = \frac{1}{n}\sum_{i=1}^n U_i^2$, $\langle\mathbf{U},\mathbf{g}\rangle_n = \frac{1}{n}\sum_{i=1}^n U_i g(Z_i)$. Notice then that, by the definition of $\gamma_n$ and $s_n$, we have

$$\gamma_n(s_n) - \gamma_n(v) = 2\langle\mathbf{W},\mathbf{f}-\mathbf{f}_{K_n}\rangle_n - \|\mathbf{f}-\mathbf{f}_{K_n}\|_n^2$$
$$+ 2\langle\mathbf{W},\mathbf{v}-\mathbf{s}\rangle_n - \|\mathbf{s}-\mathbf{v}\|_n^2 + 2\|\mathbf{f}-\mathbf{f}_{K_n}\|_n^2.$$

Let $a_n = \text{pen}(I,K_n) - \text{pen}(I_0,K_n)$. Then, by (3.10) we obtain

$$P(\hat{I}=I) \leq P\bigg(\sup_{v \in S_{I,K_n}} (2\langle\mathbf{W},\mathbf{v}-\mathbf{s}\rangle_n - \|\mathbf{s}-\mathbf{v}\|_n^2 \geq a_n/3)\bigg)$$

(3.11)
$$+ P(2\langle\mathbf{W},\mathbf{f}-\mathbf{f}_{K_n}\rangle_n - \|\mathbf{f}-\mathbf{f}_{K_n}\|_n^2 \geq a_n/3)$$

$$+ P(\|\mathbf{f}-\mathbf{f}_{K_n}\|_n^2 \geq a_n/6).$$

Recall that $\mathbf{P}_{K_n}$ is the projection matrix onto

$$L_{K_n} = \{(g(T_1),\ldots,g(T_n))' | g \in S_{K_n}\} \subset \mathbb{R}^n.$$

Define by $\mathbf{P}_{I,K_n}$ the projection matrix onto

$$L_{I,K_n} = \{(h(\underline{X}_1,T_1),\ldots,h(\underline{X}_n,T_n))' | h \in S_{I,K_n}\} \subset \mathbb{R}^n.$$

Then

$$2\langle\mathbf{W},\mathbf{v}-\mathbf{s}\rangle_n - \|\mathbf{s}-\mathbf{v}\|_n^2$$

$$= 2\|\mathbf{P}_{I,K_n}\mathbf{s} - \mathbf{v}\|_n \bigg\langle\mathbf{W},\frac{\mathbf{v}-\mathbf{P}_{I,K_n}\mathbf{s}}{\|\mathbf{P}_{I,K_n}\mathbf{s}-\mathbf{v}\|_n}\bigg\rangle_n - \|\mathbf{P}_{I,K_n}\mathbf{s}-\mathbf{v}\|_n^2$$

$$+ 2\|\mathbf{P}_{I,K_n}\mathbf{s} - \mathbf{s}\|_n \bigg\langle\mathbf{W},\frac{\mathbf{P}_{I,K_n}\mathbf{s}-\mathbf{s}}{\|\mathbf{P}_{I,K_n}\mathbf{s}-\mathbf{s}\|_n}\bigg\rangle_n - \|\mathbf{P}_{I,K_n}\mathbf{s}-\mathbf{s}\|_n^2$$

$$\leq \bigg\langle\mathbf{W},\frac{\mathbf{v}-\mathbf{P}_{I,K_n}\mathbf{s}}{\|\mathbf{P}_{I,K_n}\mathbf{s}-\mathbf{v}\|_n}\bigg\rangle_n^2 + \bigg\langle\mathbf{W},\frac{\mathbf{P}_{I,K_n}\mathbf{s}-\mathbf{s}}{\|\mathbf{P}_{I,K_n}\mathbf{s}-\mathbf{s}\|_n}\bigg\rangle_n^2$$

$$= \bigg\langle\mathbf{W}-\mathbf{P}_{I,K_n}\mathbf{W},\frac{\mathbf{v}-\mathbf{P}_{I,K_n}\mathbf{s}}{\|\mathbf{P}_{I,K_n}\mathbf{s}-\mathbf{v}\|_n}\bigg\rangle_n^2$$

$$+ \bigg\langle\mathbf{P}_{I,K_n}\mathbf{W},\frac{\mathbf{v}-\mathbf{P}_{I,K_n}\mathbf{s}}{\|\mathbf{P}_{I,K_n}\mathbf{s}-\mathbf{v}\|_n}\bigg\rangle_n^2 + \bigg\langle\mathbf{W},\frac{\mathbf{P}_{I,K_n}\mathbf{s}-\mathbf{s}}{\|\mathbf{P}_{I,K_n}\mathbf{s}-\mathbf{s}\|_n}\bigg\rangle_n^2$$

$$\leq \|\mathbf{P}_{I,K_n}\mathbf{W}\|_n^2 + \left\langle \mathbf{W}, \frac{\mathbf{P}_{I,K_n}\mathbf{s} - \mathbf{s}}{\|\mathbf{P}_{I,K_n}\mathbf{s} - \mathbf{s}\|_n} \right\rangle_n^2,$$

using $2|xy| \leq x^2 + y^2$ for the first inequality and Cauchy–Schwarz for the last one. Using an identical reasoning, we also obtain that

$$2\langle \mathbf{W}, \mathbf{f} - \mathbf{f}_{K_n}\rangle_n - \|\mathbf{f} - \mathbf{f}_{K_n}\|_n^2 \leq \|\mathbf{P}_{K_n}\mathbf{W}\|_n^2 + \left\langle \mathbf{W}, \frac{\mathbf{P}_{K_n}\mathbf{f} - \mathbf{f}}{\|\mathbf{P}_{K_n}\mathbf{f} - \mathbf{f}\|_n} \right\rangle_n^2.$$

All ratios introduced above are defined to be zero if the denominator is zero, in which case the first two probabilities in (3.11) are identically zero. Then, for the first two terms in (3.11),

(3.12)
$$P\left( \sup_{v \in S_{I,K_n}} \left( 2\langle \mathbf{W}, \mathbf{v} - \mathbf{s}\rangle_n - \|\mathbf{s} - \mathbf{v}\|_n^2 \geq \frac{a_n}{3} \right) \right)$$
$$\leq P\left( \|\mathbf{P}_{I,K_n}\mathbf{W}\|_n^2 \geq \frac{a_n}{6} \right) + P\left( \left\langle \mathbf{W}, \frac{\mathbf{P}_{I,K_n}\mathbf{s} - \mathbf{s}}{\|\mathbf{P}_{I,K_n}\mathbf{s} - \mathbf{s}\|_n} \right\rangle_n^2 \geq \frac{a_n}{6} \right).$$

Similarly,

(3.13)
$$P\left( 2\langle \mathbf{W}, \mathbf{f} - \mathbf{f}_{K_n}\rangle_n - \|\mathbf{f} - \mathbf{f}_{K_n}\|_n^2 \geq \frac{a_n}{3} \right)$$
$$\leq P\left( \|\mathbf{P}_{K_n}\mathbf{W}\|_n^2 \geq \frac{a_n}{6} \right) + P\left( \left\langle \mathbf{W}, \frac{\mathbf{P}_{K_n}\mathbf{f} - \mathbf{f}}{\|\mathbf{P}_{K_n}\mathbf{f} - \mathbf{f}\|_n} \right\rangle_n^2 \geq \frac{a_n}{6} \right).$$

We shall use Rosenthal's inequality, stated in Appendix A.4, to bound the second term in both (3.12) and (3.13). The application of this inequality, as in the proof of Theorem 3.1, page 484, of Baraud (2000), leads to

(3.14) $$P\left( \left\langle \mathbf{W}, \frac{\mathbf{P}_{I,K_n}\mathbf{s} - \mathbf{s}}{\|\mathbf{P}_{I,K_n}\mathbf{s} - \mathbf{s}\|_n} \right\rangle_n^2 \geq \frac{a_n}{6} \right) \leq C(p) E|W_1|^p n^{-p/2} \left( \frac{a_n}{6} \right)^{-p/2},$$

for any $p \geq 2$ and a constant $C(p)$ that depends only on $p$.

Recall that by the definition (3.5), $\mathrm{pen}(I, K_n) = 2(|I| + 1)rK_n \log n / n$, and so

(3.15) $$a_n = \mathrm{pen}(I, K_n) - \mathrm{pen}(I_0, K_n) \geq 2K_n \log n / n,$$

by the definition of the penalty term (3.5). Recall that, by Assumption 2.2, we have $E|W_1|^p < \infty$ for any $p > d + 4$, for some $d > 0$. Then, by (3.14) and defining $B = 3^{p/2} C(p) E|W_1|^p$, we obtain

(3.16) $$P\left( \left\langle \mathbf{W}, \frac{\mathbf{P}_{I,K_n}\mathbf{s} - \mathbf{s}}{\|\mathbf{P}_{I,K_n}\mathbf{s} - \mathbf{s}\|_n} \right\rangle_n^2 \geq \frac{a_n}{6} \right) \leq \frac{B}{(K_n \log n)^{p/2}} \to 0,$$

and, using an identical argument,

(3.17) $$P\left( \left\langle \mathbf{W}, \frac{\mathbf{P}_{K_n}\mathbf{f} - \mathbf{f}}{\|\mathbf{P}_{K_n}\mathbf{f} - \mathbf{f}\|_n} \right\rangle_n^2 \geq \frac{a_n}{6} \right) \leq \frac{B}{(K_n \log n)^{p/2}} \to 0,$$

since $K_n \asymp n^{1/2\alpha+1} \to \infty$ for all $\alpha > 0$.

We will now invoke Corollary 5.1, page 478, in Baraud (2000) to bound the first term in either (3.12) or (3.13). We discuss first (3.12). Notice that $E\|\mathbf{P}_{K_n}\mathbf{W}\|_n^2 = \sigma^2 \operatorname{tr}(\mathbf{P}_{K_n})/n = \sigma^2 r K_n/n$, using the standard properties of the projection operator. Then, by Corollary 5.1 of Baraud (2000), for any $p \geq 2$ such that $E|W_1|^p < \infty$ and for any $m > 0$,

$$(3.18)\quad P(n\|\mathbf{P}_{K_n}\mathbf{W}\|_n^2 \geq r\sigma^2 K_n + 2\sigma^2\sqrt{rK_n m} + \sigma^2 m) \leq C'(p)\vartheta_p r \frac{K_n}{m^{p/2}},$$

for some constant $C'(p) > 0$ depending only on $p$ and for $\vartheta_p = E|W_1|^p/\sigma^p$ and $\sigma^2 = EW_1^2$.

We shall use (3.18) with $m = K_n^{3/4}\log n$. Recall now (3.15) and notice that, for $n$ large enough,

$$na_n/6 \geq r\sigma^2 K_n + 2\sigma^2\sqrt{rK_n K_n^{3/4}\log n} + \sigma^2 K_n^{3/4}\log n,$$

for all $K_n \in \mathcal{K}_n$. We mention that other choices of $m$ are possible, but at the price of additional technicalities and with very little gain in terms of the overall result. Then, with $B' = rC'(4)\vartheta_4$ and using (3.18), we obtain

$$(3.19)\qquad P\left(n\|\mathbf{P}_{K_n}\mathbf{W}\|_n^2 \geq \frac{na_n}{6}\right) \leq \frac{B'}{K_n^{(3p-8)/8}(\log n)^{p/2}} \to 0,$$

since, by Assumption 2.2, $p > 4 + d > 8/3$ for any $d > 0$.

For (3.13) we first notice that $E\|\mathbf{P}_{I,K_n}\mathbf{W}\|_n^2 = \sigma^2(|I| + rK_n)/n$. Then, if we replace above $rK_n$ by $|I| + rK_n$, we also obtain

$$(3.20)\qquad P\left(n\|\mathbf{P}_{I,K_n}\mathbf{W}\|_n^2 \geq \frac{na_n}{6}\right) \leq \frac{B''}{K_n^{(3p-8)/8}(\log n)^{p/2}} \to 0,$$

for an appropriately modified constant $B''$.

We bound now the last term in (3.11). Recall the approximation error bound of (2.7). By Markov's inequality and (3.15) we have, with $C = 3C^2(\alpha)L^2 h_1$,

$$(3.21)\quad \begin{aligned} P\left(\|\mathbf{f} - \mathbf{f}_{K_n}\|_n^2 \geq \frac{a_n}{6}\right) &\leq \frac{6E\|f - f_{K_n}\|_{\mu_T}^2}{a_n} \\ &\leq \frac{C}{K_n^{2\alpha}} \times \frac{n}{K_n \log n} \leq \frac{C}{\log n} \to 0, \end{aligned}$$

by the choice of $K_n \asymp n^{1/2\alpha+1}$.

Notice now that the number of terms in (3.8) is bounded by $A_1$, where $A_1 > 0$, independent of $n$, is the number of models the linear part of which includes the $I_0$ variables. Then, by (3.8), (3.11), (3.16), (3.17) and (3.19)–(3.21), we obtain $P(I_0 \subsetneq \hat{I}) \to 0$, which is the desired result.

(b) The proof is almost identical for this case, in which we now re-denote $K_{n,a} \asymp n^{1/2a+2}$ by $K_n$ and use the penalty term (3.6) instead of (3.5). Note that (3.16), (3.17), (3.19) and (3.20) only require that $K_n \to \infty$, and thus they hold independently of $\alpha$ or $a$. The only difference is in (3.21), which now becomes

$$P\left(\|\mathbf{f} - \mathbf{f}_{K_n}\|_n^2 \geq \frac{a_n}{6}\right) \leq \frac{6E\|f - f_{K_n}\|_{\mu_T}^2}{a_n}$$

$$\leq \frac{C}{K_n^{2\alpha}} \times \frac{n}{K_n \log n} \leq \frac{C}{K_n^{2a+1}} \times \frac{n}{K_n \log n} \quad \text{(uniformly over } \alpha \geq 1/2 + a\text{)}$$

$$\leq \frac{C}{\log n} \to 0 \quad \text{for } K_{n,a} \asymp n^{1/2a+2}.$$

Then the concluding argument is identical to the one above.

2. $P(I_0 \not\subset \hat{I}) \to 0$. The proof is the same in both cases. Let $c = \inf_{j \in I_0} |\beta_j|$ and notice that, by the definition of $I_0$, we have $c > 0$. Consequently, $|\hat{\beta}_{j_n} - \beta_{j_n}| = |\beta_{j_n}|$, for all $j_n \in I_0 \setminus \hat{I}$, and

$$P(I_0 \not\subset \hat{I}) \leq P(j_n' \notin \hat{I}, \text{ for some } j_n' \in I_0)$$

$$\leq P(|\hat{\beta}_{j_n'} - \beta_{j_n'}| = |\beta_{j_n'}|) \leq P(|\hat{\beta}_{j_n'} - \beta_{j_n'}| \geq c) \to 0,$$

by the component-wise consistency of $\hat{\beta}_{\hat{I}}$ implied by 3 of Corollary 2.1. This completes the proof of this theorem. □

## 4. Asymptotic normality.

4.1. *Asymptotic normality of $\hat{\beta}_{I_0}$.* In this section we assume that $I_0 \subseteq \{1, \ldots, q\}$ is known. Then, with the notation of Section 2.2, (1.1) becomes

(4.1) $$Y_i = \beta_{I_0}' X_{I_0, i} + f(T_i) + W_i.$$

Here $X_{I_0}$ denotes the vector of covariates corresponding to the index set $I_0$. Let $|I_0| = q_0$ for some known $q_0 \leq q$. In order to emphasize the dimension of the parametric part, we re-denote $\beta_{I_0}$ by $\beta_{q_0}$ and its estimator by $\hat{\beta}_{q_0}$. We first consider estimators for $\beta_{q_0}$ and $f$ within the family of approximating spaces $\{S_{I_0, K_n}\}_{K_n \in \mathcal{K}_n}$. Recall that this corresponds to our Case 1 defined in the Introduction. We show in the following theorem that estimating $f$ adaptively preserves the asymptotic normality of $\hat{\beta}_{q_0}$.

Let $\Sigma_{q_0} = (\sigma_{kj})_{q_0 \times q_0}$, with $\sigma_{kj} = \text{Cov}(X_k, X_j) - \text{Cov}(\theta_k(T), \theta_j(T))$, for $(j, k) \in I_0 \times I_0$, and let $\text{Var}(W) = \sigma^2$.

THEOREM 4.1. *Let $\hat{\beta}_{q_0}$ be the estimator of $\beta_{q_0}$ based on the approximating spaces $\{S_{I_0, K_n}\}_{K_n \in \mathcal{K}_n}$. Under Assumptions 2.1–2.5, if $f \in \mathcal{D}_{\alpha, 2}(L)$, then*

$$\sqrt{n}(\hat{\beta}_{q_0} - \beta_{q_0}) \xrightarrow{d} N_{q_0}(0, \sigma^2 \Sigma_{q_0}^{-1}).$$

Theorem 4.1 holds uniformly over $\alpha \in (1/2, r)$ and $\gamma \in (b/4, r)$, $b \geq 3$.

Corollary 4.1 is an immediate consequence of this theorem. Recall the definitions $K_{n,\alpha} = n^{1/2\alpha+1}$ and $K_{n,a} = n^{1/2a+2}$.

COROLLARY 4.1. *Let $\hat{\beta}_{q_0}$ be the estimator of $\beta_{q_0}$ based on either $S_{I_0, K_{n,\alpha}}$ or $S_{I_0, K_{n,a}}$. Under Assumptions* 2.1–2.5, *if $f \in \mathcal{D}_{\alpha,2}(L)$, then*

$$\sqrt{n}(\hat{\beta}_{q_0} - \beta_{q_0}) \xrightarrow{d} N_{q_0}(0, \sigma^2 \Sigma_{q_0}^{-1}).$$

Corollary 4.1 holds uniformly over $\gamma \in (b/4, r)$, $b \geq 3$. Also, it holds for any *fixed* $\alpha \in (1/2, r)$ if $S_{I_0, K_{n,\alpha}}$ is used, and uniformly over $\alpha \in (1/2 + a, r)$ for $S_{I_0, K_{n,a}}$. We note that the estimators $\hat{\beta}_{q_0}$ are different in each situation; we have used the same notation for brevity, since we are only interested in their limiting distribution. The proofs are given in Appendix A.3.

4.2. *Cases* 2 *and* 3: *asymptotic normality of $\hat{\beta}_{\hat{I}}$.* Throughout this section we regard $\hat{\beta}_{\hat{I}}$ and $\hat{\beta}_{I_0}$ as vectors in $\mathbb{R}^q$, by adding 0's in the necessary positions. Also, we re-denote $\beta \in \mathbb{R}^q$ by $\beta_{I_0} \in \mathbb{R}^q$ to emphasize that the only nonzero components of $\beta$ correspond to the index set $I_0$. Let $\mathbf{V}^0 = \sigma^2 \Sigma_{q_0}^{-1}$. Let $\mathbf{V} = (\mathbf{V}_{ij})_{q \times q}$, where $\mathbf{V}_{ij} = \mathbf{V}_{ij}^0$ for $(i, j) \in I_0 \times I_0$ and zero otherwise.

THEOREM 4.2. *If $f \in \mathcal{D}_{\alpha,2}(L)$ and Assumptions* 2.1–2.5 *hold, then we have the following results:*

*Case* 2. *If* $\text{pen}(I, K_{n,\alpha}) = 2(|I| + 1) K_{n,\alpha} \log n / n$, *for some given* $1/2 < \alpha < r$, *then* $\sqrt{n}(\hat{\beta}_{\hat{I}} - \beta_{I_0}) \xrightarrow{d} N_q(0, \mathbf{V})$.

*Case* 3. *If* $\text{pen}(I, K_n) = 2(|I|+1) K_{n,a} \log n / n$, *then* $\sqrt{n}(\hat{\beta}_{\hat{I}} - \beta_{I_0}) \xrightarrow{d} N_q(0, \mathbf{V})$, *uniformly over* $\alpha \in [1/2 + a, r)$ *for* $0 < a < r - 1/2$.

*In both cases, the limiting distribution has all its mass concentrated on the space generated by the $I_0$ covariates.*

PROOF. We prove that, in both cases, for any $c \in \mathbb{R}^q$,

(4.2) $$c'\sqrt{n}(\hat{\beta}_{\hat{I}} - \beta_{I_0}) \xrightarrow{d} N(0, c'\mathbf{V}c) \qquad \text{as } n \to \infty,$$

which leads to the desired result. For any $b \in \mathbb{R}$, $c \in \mathbb{R}^q$ we have

$$P(c'(\sqrt{n}(\hat{\beta}_{\hat{I}} - \beta_{I_0})) \leq b)$$

(4.3) $$= P(c'(\sqrt{n}(\hat{\beta}_{\hat{I}} - \beta_{I_0})) \leq b, \hat{I} = I_0)$$
$$+ P(c'(\sqrt{n}(\hat{\beta}_{\hat{I}} - \beta_{I_0})) \leq b, \hat{I} \neq I_0)$$
$$= P(c'(\sqrt{n}(\hat{\beta}_{I_0} - \beta_{I_0})) \leq b) - P(c'(\sqrt{n}(\hat{\beta}_{I_0} - \beta_{I_0})) \leq b, \hat{I} \neq I_0)$$
$$+ P(c'(\sqrt{n}(\hat{\beta}_{\hat{I}} - \beta_{I_0})) \leq b, \hat{I} \neq I_0).$$

Since, by definition, $\hat{\beta}_{I_0} \in \mathbb{R}^q$ has nonzero elements only in the positions corresponding to $I_0$, then

(4.4) $$c'(\sqrt{n}(\hat{\beta}_{I_0} - \beta_{I_0})) = c'_0(\sqrt{n}(\hat{\beta}_{q_0} - \beta_{q_0})),$$

where $c_0$, $\hat{\beta}_{q_0}$ and $\beta_{q_0}$ are obtained from $c$, $\hat{\beta}_{I_0}$ and $\beta_{I_0}$, respectively, by deleting the coordinates corresponding to zeros in $\beta_{I_0}$. Now, we have that

$$P(c'(\sqrt{n}(\hat{\beta}_{\hat{I}} - \beta_{I_0})) \leq b, \hat{I} \neq I_0) \leq P(\hat{I} \neq I_0)$$

and that

$$P(c'(\sqrt{n}(\hat{\beta}_{I_0} - \beta_{I_0})) \leq b, \hat{I} \neq I_0) \leq P(\hat{I} \neq I_0).$$

By Theorem 3.1, $P(\hat{I} \neq I_0) \to 0$ in both cases. Thus, from (4.3) and (4.4) we obtain

(4.5) $$\lim_{n \to \infty} P(\sqrt{n}(c'(\hat{\beta}_{\hat{I}} - \beta_{I_0})) \leq b) = \lim_{n \to \infty} P(\sqrt{n}(c'_0(\hat{\beta}_{q_0} - \beta_{q_0})) \leq b).$$

In both cases, the right-hand side in (4.5) converges to $N(0, c'_0 \sigma^2 \Sigma_{q_0}^{-1} c_0)$, by Corollary 4.1. Lemma A.5, in Appendix A.4, shows that the only symmetric and semipositive definite matrix $\mathbf{V}$ such that

$$c'_0 \mathbf{V}^0 c_0 = c' \mathbf{V} c \qquad \text{for any } c \in \mathbb{R}^q,$$

is given by $\mathbf{V}_{ij} = \mathbf{V}^0_{ij}$ for $(i, j) \in I_0 \times I_0$ and zero otherwise. Then (4.2) holds and the proof of the theorem is complete. □

We discuss below the limiting covariance matrix given by Theorems 4.1 and 4.2, under the assumption that the error distribution is Gaussian.

First note that the information bound for a regular estimator of $\beta$ in (4.1) is $\sigma^2 \Sigma_{q_0}^{-1}$, for some known $q_0 \leq q$; see, for example, Example 5, page 110, in Bickel, Klaassen, Ritov and Wellner (1993). Then, Theorem 5.1 shows that, for *known* $I_0$, $\hat{\beta}_{q_0}$ is asymptotically efficient.

However, if $I_0$ itself is regarded as a parameter, as in our Cases 2 and 3, the classical information bound theory no longer applies and we need to resort to other means to assess the performance of our estimators. We thus verify whether our method, which is not based on a priori knowledge of $I_0$, leads to estimators with the same limit behavior as of those constructed knowing $I_0$. By Theorem 4.2 and the continuous mapping theorem, we obtain

(4.6) $$\sqrt{n}(\hat{\beta}_{q_0} - \beta_{q_0}) \xrightarrow{d} N_{q_0}(0, \sigma^2 \Sigma_{q_0}^{-1})$$

and $\sqrt{n}(\hat{\beta}_{q_1} - \beta_{q_1}) \xrightarrow{d} 0$, where $\beta_{q_1}$ denotes the vector of zero coefficients in $\beta$. Then, indeed, $\hat{\beta}_{q_0}$ achieves the information bound for $\beta_{q_0}$ in (4.1). Notice now that if $I_0$ is known prior to estimation, then one can set $\hat{\beta}_{q_1}$ to zero, whereas our method may estimate $\beta_{q_1}$ by nonzero sequences. However, they converge to zero at an $n^{-1/2}$ rate.

For Cases 2 and 3 a much simpler method of estimation in (1.1) is to fit the model with all covariates included. This will reduce the computing time considerably, since we will only fit one model. We denote by $\tilde{\beta} \in \mathbb{R}^q$ the corresponding estimator; note that this estimator is different in the two cases, but here we are only interested in its limiting distribution, so we use the same notation. This simpler procedure of estimation is also independent of $I_0$ and, as above, we study the performance of $\tilde{\beta}$ by investigating its asymptotic properties under the assumption that $I_0$ is known. Using the same reasoning as in Theorem A.1 and, with $\beta_{I_0}$ denoting, as above, a vector in $\mathbb{R}^q$ having nonzero components only in the $I_0$ positions, one can show that $\sqrt{n}(\tilde{\beta} - \beta_{I_0}) \xrightarrow{d} N_q(0, \sigma^2 \Sigma_q^{-1})$. Let $\mathbf{I} = \Sigma_q/\sigma^2$. Consider the partition of $\mathbf{I}$ in blocks $\mathbf{I}_{11}, \mathbf{I}_{12}, \mathbf{I}_{21}, \mathbf{I}_{22}$, where $\mathbf{I}_{11}$ and $\mathbf{I}_{22}$ have dimensions $q_0 \times q_0$ and $q_1 \times q_1$, respectively. Then, as in Bickel, Klaassen, Ritov and Wellner [(1993), page 28],

$$(4.7) \qquad \sqrt{n}(\tilde{\beta}_{q_0} - \beta_{q_0}) \xrightarrow{d} N_{q_0}(0, \mathbf{I}_{11.2}^{-1})$$

and $\sqrt{n}(\tilde{\beta}_{q_1} - \beta_{q_1}) \to N_{q_1}(0, \mathbf{I}_{22.1}^{-1})$, where $\mathbf{I}_{11.2} = \mathbf{I}_{11} - \mathbf{I}_{12}\mathbf{I}_{22}^{-1}\mathbf{I}_{21}$ and $\mathbf{I}_{22.1} = \mathbf{I}_{22} - \mathbf{I}_{21}\mathbf{I}_{11}^{-1}\mathbf{I}_{12}$. Note that the limiting distribution in (4.7) coincides with the one in (4.6) only if $\mathbf{I}_{12} = 0$. Thus, although this procedure might be more appealing from a computational point of view, it leads to estimators with higher variance, if the true model corresponds to a proper subset of the full collection of the X covariates.

**5. Conclusions.** This article studies simultaneous estimation of $\beta$, $I_0$ and $f$. We showed that one can consistently estimate $I_0$ and obtain asymptotically normal estimators for the selected estimator of $\beta$. The construction of the approximating spaces used for estimation parallels the one used in parametric or nonparametric model selection problems, but the penalty term needs to be adjusted for the semiparametric case. We summarize our findings in each of the cases under consideration.

- Case 2: $I_0$ unknown and $\alpha > 1/2$ known. O1 and O2 hold. For the approximating spaces $\{S_{I,K_{n,\alpha}}\}_{I \in \mathcal{I}}$ and the penalty term $\text{pen}(I, K_{n,\alpha}) = 2(|I| + 1)rK_{n,\alpha}\log n/n$, with $K_{n,\alpha} \asymp n^{1/2\alpha+1}$, we showed that $P(\hat{I} = I_0) \to 0$ and that $\sqrt{n}(\hat{\beta}_{\hat{I}} - \beta_{I_0}) \xrightarrow{d} N_q(0, \mathbf{V})$ *for a specified* $\alpha$. The rate of convergence $r_n$ of $\hat{f}$ is of order $O(\log n/n^{2\alpha/2\alpha+1})$, which is the minimax optimal rate for a given $\alpha$, up to a $\log n$ factor.
- Case 3: $I_0$ and $\alpha > 1/2$ unknown. O1 and O2 hold. For the approximating spaces $\{S_{I,K_{n,a}}\}_{I \in \mathcal{I}}$, with $K_{n,a} \asymp n^{1/2a+2}$, $a > 0$ arbitrarily close to zero, and for the penalty term $\text{pen}(I, K_{n,a}) = 2(|I|+1)rK_{n,a}\log n/n$, we showed that $P(\hat{I} = I_0) \to 0$ and that $\sqrt{n}(\hat{\beta}_{\hat{I}} - \beta_{I_0}) \xrightarrow{d} N_q(0, \mathbf{V})$ *uniformly over* $\alpha > 1/2$. The rate $r_n$ of $\hat{f}$ is of order $O(\log n/n^{(2a+1)/(2a+2)})$.

We also note that in Case 1: $I_0$ known, $\alpha$ unknown, for the approximating spaces $\{S_{I_0,K_n}\}_{K_n \in \mathcal{K}_n}$ and the penalty term $\text{pen}(I, K_n) = 2(|I| + 1)rK_n \log n/n$, we have shown that $\sqrt{n}(\hat{\beta}_{I_0} - \beta_{I_0}) \xrightarrow{d} N_{q_0}(0, \sigma^2 \Sigma_{q_0}^{-1})$ for $q_0 = |I_0|$. Also, $r_n = O((\log n/n)^{2\alpha/2\alpha+1})$ for any $\alpha \in (1/2, r)$. Thus, $\hat{f}$ is minimax adaptive, up to a $\log n$ factor.

## APPENDIX

**A.1. The unicity of the least squares estimators.** In this section we give the proof of the asymptotic unicity of our estimators.

LEMMA A.1. *Under Assumptions* 2.1 *and* 2.3, $\mathbf{Z}'_{K_n} \mathbf{Z}_{K_n}$ *is invertible for any* $K_n \in \mathcal{K}_n$, *except for an event whose probability tends to zero as* $n \to \infty$.

PROOF. Notice that, for any $K_n \in \mathcal{K}_n$,

$$\mathbf{Z}'_{K_n} \mathbf{Z}_{K_n} = \left( \sum_{i=1}^n \phi_j(T_i) \phi'_{j'}(T_i) \right)_{1 \le j, j' \le rK_n}.$$

Then $\mathbf{Z}'_{K_n} \mathbf{Z}_{K_n}$ is invertible if and only if $\mathbf{\Phi} = n^{-1} \mathbf{Z}'_{K_n} \mathbf{Z}_{K_n}$ is invertible. Let $\{\zeta_j\}_{j=1}^{rK_n}$ be the eigenvalues of $\mathbf{\Phi}$, and let $(1/\zeta)_{\max} = \sup\{1/\zeta_1, \ldots, 1/\zeta_{rK_n}\}$. Thus, $\mathbf{\Phi}$ is invertible on the set where $(1/\zeta)_{\max} \le \rho_0 < \infty$, for some $\rho_0 > 0$. By the proof of Lemma 3.1, page 492, in Baraud (2000), since $(\phi_j)_{j=1}^{rK_n}$ are orthogonal in $L_2(\lambda_T, [0,1])$, we have

$$\left(\frac{1}{\zeta}\right)_{\max} = \sup_{u \in S_{K_n} \setminus \{0\}} \frac{\|u\|_{\lambda_T}^2}{\|u\|_n^2}.$$

Notice that Assumption 2.1 implies that the density of $T$ with respect to $\lambda_T$ is bounded above and below by $Lh_1$ and $Lh_0$, respectively, where $0 < L < \infty$ is the Lebesgue measure of the compact set $K$. Also, under Assumptions 2.1 and 2.3, Proposition A.2, adapted to the case $q = 0$, ensures that for some constant $D_1 > 0$ condition $H_{\text{Con}} : \|u\|_\infty \le D_1 \sqrt{rN_n} \|u\|_{\lambda_T}$ for all $u \in S_{N_n}$ of Baraud (2002) holds. Then, by Lemma 6.2, page 21, and the proof of Proposition 5.2, page 24, in Baraud (2002), for all $\rho_0 > L^{-1} h_0^{-1}$ and for a constant $D > 0$ depending on $L, h_0, h_1$ and $\rho_0$, we obtain

$$P\left(\left(\frac{1}{\zeta}\right)_{\max} > \rho_0\right) = P\left(\sup_{u \in S_{K_n} \setminus \{0\}} \frac{\|u\|_{\lambda_T}^2}{\|u\|_n^2} > \rho_0\right)$$

$$\le P\left(\sup_{u \in S_{N_n} \setminus \{0\}} \frac{\|u\|_{\lambda_T}^2}{\|u\|_n^2} > \rho_0\right)$$

$$\leq r^2 N_n^2 \exp\left(-\frac{Dn}{D_1^2 N_n^2}\right) \to 0,$$

where the first inequality holds because the approximating spaces are nested and the convergence to zero holds since $N_n \asymp (n/\log n)^{1/2}$, by construction. This concludes the proof of this proposition. □

REMARK A.1. Since $K_{n,\alpha} \asymp n^{1/2\alpha+1}$ and $K_{n,a} = n^{1/2a+2}$ belong to $\mathcal{K}_n$, by construction (see page 903) the above result implies that $\mathbf{Z}'_{K_{n,\alpha}} \mathbf{Z}_{K_{n,\alpha}}$ and $\mathbf{Z}'_{K_{n,a}} \mathbf{Z}_{K_{n,a}}$ are also invertible.

LEMMA A.2. *Under Assumptions* 2.1–2.5, *for any* $I \subseteq \{1,\ldots,q\}$, $\mathbf{X}'_I(\mathbf{Id} - \mathbf{P}_{\hat{K}_n})\mathbf{X}_I$ *and* $\mathbf{X}'_I(\mathbf{Id} - \mathbf{P}_{K_{n,\alpha}})\mathbf{X}_I$, $\alpha \in (1/2, r)$, *are invertible except for an event whose probability tends to zero as* $n \to \infty$.

The proof of this lemma is based on Proposition A.1, which in turn requires the proof of Lemma A.3. We use the following notation here and in the sequel.

NOTATION A.1. Let $\Sigma_{|I|} = (\sigma_{kj})_{|I| \times |I|}$, with $\sigma_{kj} = \text{Cov}(X_k, X_j) - \text{Cov}(\theta_k(T), \theta_j(T))$, for $(j,k) \in I \times I$, $I \in \mathcal{I}$. Let $\boldsymbol{\theta}_j = (\theta_j(T_1),\ldots,\theta_j(T_n))'$, $\varepsilon_{ij} = X_{ij} - \theta_j(T_i)$, $\boldsymbol{\varepsilon}_j = (\varepsilon_{1j},\ldots,\varepsilon_{nj})'$, for $j \in I$ and $1 \leq i \leq n$. Also, let $\boldsymbol{\theta}$ and $\boldsymbol{\varepsilon}$ be the $n \times |I|$ matrices having columns $\boldsymbol{\theta}_j$ and $\boldsymbol{\varepsilon}_j$, respectively. Let $I_q = \{1,\ldots,q\}$. Recall that for any $\mathbf{U} = (U_1,\ldots,U_n)'$ we denote $\|\mathbf{U}\|_n^2 = \frac{1}{n}\sum_{i=1}^n U_i^2$.

LEMMA A.3. *Under Assumptions* 2.1–2.5,

(A.1) $$\|(\mathbf{Id} - \mathbf{P}_{\hat{K}_n})\boldsymbol{\theta}_j\|_n = O_P((\log n/n)^{1/4}),$$

*for any* $j \in I_q$.

PROOF. Let $b_n = (\log n/n)^{1/2}$. First notice that

(A.2) $$P(\|(\mathbf{Id} - \mathbf{P}_{\hat{K}_n})\boldsymbol{\theta}_j\|_n^2 \geq b_n) \leq \sum_{k_n=A_n}^{J_n} P(\|(\mathbf{Id} - \mathbf{P}_{2^{k_n}})\boldsymbol{\theta}_j\|_n^2 \geq b_n),$$

with $K_n = 2^{k_n}$. Next, observe that $\mathbf{P}_{K_n}\boldsymbol{\theta}_j$ is the projection of $\boldsymbol{\theta}_j$ onto $L_{K_n} = \{(g(T_1),\ldots,g(T_n))'; g \in S_{K_n}\}$. For every $K_n \in \mathcal{K}_n$ let $\theta_{j,K_n}$ be the $L_2(\mu_T)$ projection of $\theta_j$ onto $S_{K_n}$. Then

$$\|(\mathbf{Id} - \mathbf{P}_{K_n})\boldsymbol{\theta}_j\|_n \leq \|\boldsymbol{\theta}_j - \boldsymbol{\theta}_{j,K_n}\|_n.$$

By Assumption 2.5 $\theta_j \in \mathcal{D}_{\gamma,2}(A)$. Thus, by (2.7) $\|\theta_j - \theta_{j,K_n}\|_{\mu_T}^2 \leq h_1 C(\gamma, A) \times K_n^{-2\gamma}$. Then, by Markov's inequality and recalling, by (2.1), that $A_n =$

$[\log_2(n/\log n)^{1/b}]$, $b \geq 3$, and denoting $J_n = [\log_2(n/\log n)^{1/2}]$, we obtain

$$P(\|(\mathbf{Id} - \mathbf{P}_{\hat{K}_n})\boldsymbol{\theta}_j\|_n^2 \geq b_n) \leq \sum_{k_n=A_n}^{J_n} P(\|\boldsymbol{\theta}_j - \boldsymbol{\theta}_{j,K_n}\|_n^2 \geq b_n)$$

$$\leq \sum_{k_n=A_n}^{J_n} h_1 C(\gamma, A) \frac{1}{K_n^{2\gamma}} \frac{n^{1/2}}{(\log n)^{1/2}}$$

$$\leq h_1 C(\gamma, A) \log_2 n \frac{(\log n)^{2\gamma/b - 1/2}}{n^{2\gamma/b - 1/2}}$$

$$\to 0,$$

for $\gamma > b/4$, which concludes the proof of this result. $\square$

REMARK A.2. The proof of the above result also implies that (A.1) holds with $\hat{K}_n$ replaced by any $K_n \in \mathcal{K}_n$. Thus, in particular, it holds for $K_{n,\alpha} \asymp n^{1/2\alpha+1}$ and $K_{n,a} = n^{1/2a+2}$.

PROPOSITION A.1. *Under Assumptions* 2.1–2.5 *we have*

(A.3) $$\mathbf{X}_I'(\mathbf{Id} - \mathbf{P}_{\hat{K}_n})\mathbf{X}_I/n \xrightarrow{P} \Sigma_{|I|}.$$

PROOF. With the Notation A.1, each row in $\mathbf{X}_I'$ can be written as $\boldsymbol{\theta}_j' + \boldsymbol{\varepsilon}_j'$. Then, for every $j, k \in I$, we have

(A.4)
$$\begin{aligned}
(\mathbf{X}_I'(\mathbf{Id} - \mathbf{P}_{\hat{K}_n})\mathbf{X}_I)_{jk} \\
= (\boldsymbol{\varepsilon}_j + \boldsymbol{\theta}_j)'(\mathbf{Id} - \mathbf{P}_{\hat{K}_n})(\boldsymbol{\varepsilon}_k + \boldsymbol{\theta}_k) \\
= \boldsymbol{\varepsilon}_j'(\mathbf{Id} - \mathbf{P}_{\hat{K}_n})\boldsymbol{\varepsilon}_k + \boldsymbol{\varepsilon}_j'(\mathbf{Id} - \mathbf{P}_{\hat{K}_n})\boldsymbol{\theta}_k \\
+ \boldsymbol{\theta}_j'(\mathbf{Id} - \mathbf{P}_{\hat{K}_n})\boldsymbol{\varepsilon}_k + \boldsymbol{\theta}_j'(\mathbf{Id} - \mathbf{P}_{\hat{K}_n})\boldsymbol{\theta}_k.
\end{aligned}$$

Notice that $\boldsymbol{\varepsilon}_j'\boldsymbol{\varepsilon}_k/n \xrightarrow{P} \sigma_{jk}$, by the law of large numbers and the definition of $\sigma_{jk}$. Also, as in the proof of Lemma 5 of Chen (1988),

(A.5) $$P(|n^{-1}(\boldsymbol{\varepsilon}_i'\mathbf{P}_{\hat{K}_n}\boldsymbol{\varepsilon}_j)| \geq c) \leq rc^{-1}n^{-1}E(\hat{K}_n) \to 0$$

for any $c > 0$, since $E(\hat{K}_n) \leq N_n \asymp (n/\log n)^{1/2}$. Notice now that, by the Cauchy–Schwarz inequality,

$$n^{-1}|\boldsymbol{\varepsilon}_j'(\mathbf{Id} - \mathbf{P}_{\hat{K}_n})\boldsymbol{\theta}_k| \leq \|\boldsymbol{\varepsilon}_j\|_n \|(\mathbf{Id} - \mathbf{P}_{\hat{K}_n})\boldsymbol{\theta}_k\|_n \xrightarrow{P} 0,$$

since $\|\boldsymbol{\varepsilon}_j\|_n \to_{\text{a.s}} \sigma_{jj}$, which is finite, and $\|(\mathbf{Id} - \mathbf{P}_{\hat{K}_n})\boldsymbol{\theta}_k\|_n \xrightarrow{P} 0$ by Lemma A.3. By symmetry, we also have $n^{-1}\boldsymbol{\theta}_j'(\mathbf{Id} - \mathbf{P}_{\hat{K}_n})\boldsymbol{\varepsilon}_k \xrightarrow{P} 0$. For the last term

in (A.4), by the Cauchy–Schwarz inequality and by Lemma A.3, and since $\mathbf{Id} - \mathbf{P}_{\hat{K}_n}$ is idempotent, we have

$$
\begin{aligned}
n^{-1}&|\boldsymbol{\theta}'_j(\mathbf{Id} - \mathbf{P}_{\hat{K}_n})\boldsymbol{\theta}_k| \\
&= n^{-1}|((\mathbf{Id} - \mathbf{P}_{\hat{K}_n})\boldsymbol{\theta}_j)'((\mathbf{Id} - \mathbf{P}_{\hat{K}_n})\boldsymbol{\theta}_k)| \\
&\leq \|(\mathbf{Id} - \mathbf{P}_{\hat{K}_n})\boldsymbol{\theta}_j\|_n \|(\mathbf{Id} - \mathbf{P}_{\hat{K}_n})\boldsymbol{\theta}_k\|_n \\
&\xrightarrow{P} 0.
\end{aligned}
\qquad (A.6)
$$
□

REMARK A.3. By Remark A.2 and the proof above, we can conclude that (A.3) holds with $\hat{K}_n$ replaced by any $K_n \in \mathcal{K}_n$. Thus, as before, it holds for $K_{n,\alpha} \asymp n^{1/2\alpha+1}$ and $K_{n,a} = n^{1/2a+2}$.

PROOF OF LEMMA A.2. Let $\mathbf{P}$ denote any of the two projection matrices. We show that, for any $c_I \in \mathbb{R}^{|I|} \setminus \{0\}$ the sequence $c'_I \mathbf{X}'_I(\mathbf{Id} - \mathbf{P})\mathbf{X}_I c_I$ tends in probability to $\infty$, which implies that the corresponding matrix $\mathbf{P}$ is positive definite, hence invertible, except for a set of probability tending to zero.

By Proposition A.1 and Remark A.3, we have that $\mathbf{X}'_I(\mathbf{Id} - \mathbf{P})\mathbf{X}_I/n \xrightarrow{P} \Sigma_{|I|}$ for any of the two projection matrices. Then, by the continuous mapping theorem, for any $c_I \in \mathbb{R}^I \setminus \{0\}$ we have that $c'_I \mathbf{X}'_I(\mathbf{Id} - \mathbf{P})\mathbf{X}_I c_I/n \xrightarrow{P} c'_I \Sigma_{|I|} c_I$. If we denote $\Sigma = \Sigma_{|I_q|}$, then, by Assumption 2.3, for any nonzero $l \in \mathbb{R}^q$,

$$
\begin{aligned}
l'\Sigma l &= \operatorname{Var}(l'(\underline{X} - E(\underline{X}|T))) \\
&= \int_{[0,1]} \operatorname{Var}(l'\underline{X}|t)\,d\mu_T(t) > 0.
\end{aligned}
$$

Thus $\Sigma$ is positive definite, and so is $\Sigma_{|I|}$, for any $I \in \mathcal{I}$. Hence $c'_I \mathbf{X}'_I(\mathbf{Id} - \mathbf{P})\mathbf{X}_I c_I \xrightarrow{P} \infty$, which completes the proof of this lemma. □

**A.2. The consistency of the penalized least squares estimators and rates of convergence.** We first establish Theorem 2.1. This theorem is a direct consequence of Theorem 2.1 of Baraud (2002). The next proposition verifies that its condition ($H_{\text{Con}}$) holds.

PROPOSITION A.2. *If Assumptions* 2.1 *and* 2.3 *hold, then there exists a constant $K \geq 1$ such that, for any $g \in S_{I_q,N_n}$,*

$$\|g\|_\infty \leq K\sqrt{q + rN_n}\|g\|_\lambda,$$

*where $\|g\|_\infty = \sup_{z \in \mathcal{K}} |g(z)|$. Also note that, by Assumption* 2.1 *and by construction, $g \in L_2(\mathcal{K} \times [0,1], \lambda)$. Here we recall that $S_{I_q,N_n}$ is the largest of the approximating spaces introduced in Section* 2.2, *with $I_q = \{1,\ldots,q\}$ and $N_n$ given in* (2.1).

PROOF. Recall that $\{\phi_j\}_{j=1}^{rN_n}$ is an orthonormal basis in $L_2(\lambda_T, [0,1])$ for the linear space $S_{N_n}$ defined in Section 2.2. Then we can write any $g \in S_{I_q,N_n}$ as $g(\underline{x},t) = \sum_{j=1}^{q} a_j x_j + \sum_{j=1}^{rN_n} b_j \phi_j(t)$. By Lemma 1, page 337, in Birgé and Massart (1998), we have that $\max_{t \in [0,1]} \sum_{j=1}^{rN_n} \phi_j^2(t) \leq (2r-1)^2 N_n$. Also, by Assumption 2.1, there exists an $M > 0$ such that $|\underline{X}|_2 \leq M$ with probability 1. Hence, applying the Cauchy–Schwarz inequality, we obtain

$$\|g\|_\infty^2 \leq 2 \sum_{j=1}^{q} a_j^2 \max_{\underline{X} \in \mathcal{K}} \sum_{j=1}^{q} x_j^2 + 2 \sum_{j=1}^{rN_n} b_j^2 \max_{t \in [0,1]} \sum_{j=1}^{rN_n} \phi_j^2(t)$$

(A.7)

$$\leq 2(M^2 + (2r-1)^2)(q + rN_n)\left\{ \sum_{j=1}^{q} a_j^2 + \sum_{j=1}^{rN_n} b_j^2 \right\}.$$

Next, we show that there exists a $K_1 \geq 1$ such that

(A.8) $$\sum_{j=1}^{q} a_j^2 + \sum_{j=1}^{rN_n} b_j^2 \leq K_1 \|g\|_\lambda^2.$$

Recall that $\theta_j(t) = E(X_j | T = t)$, $j = 1, \ldots, q$. Since $E\{(X_j - \theta_j(T))\phi_k(T)\} = 0$ for all $j$ and $k$, we obtain

$$\|g\|_\mu^2 = E\left\{ \sum_{j=1}^{q} a_j X_j + \sum_{j=1}^{rN_n} b_j \phi_j(T) \right\}^2$$

(A.9)

$$= E\left\{ \sum_{j=1}^{q} a_j (X_j - \theta_j(T)) \right\}^2 + E\left\{ \sum_{j=1}^{q} a_j \theta_j(T) + \sum_{j=1}^{rN_n} b_j \phi_j(T) \right\}^2.$$

Recall that $\Sigma = (\sigma_{ij})_{q \times q}$, with $\sigma_{ij} = \text{Cov}(X_i, X_j) - \text{Cov}(\theta_i(T), \theta_j(T))$, is positive definite and so its smallest eigenvalue $\lambda_{\min} > 0$, by Assumption 2.3. Also, recall that, under Assumption 2.1, $\mu$ has on its support a density with respect to $\lambda$ that is bounded below by $h_0 > 0$ and above by $h_1 < \infty$. Then, from (A.9) and under Assumption 2.1, with $\underline{a} = (a_1, \ldots, a_q)$, it follows that

$$h_1 \|g\|_\lambda^2 \geq \|g\|_\mu^2 \geq E\left\{ \sum_{j=1}^{q} a_j (X_j - \theta_j(T)) \right\}^2 = \underline{a}' \Sigma \underline{a} \geq \lambda_{\min} \sum_{j=1}^{q} a_j^2.$$

Let $\lambda_X$ denote the restriction of $\lambda$ to the compact set $K$. Then

$$\sum_{j=1}^{rN_n} b_j^2 = \left\| \sum_{j=1}^{rN_n} b_j \phi_j \right\|_{\lambda_T}^2 \quad \text{[since } \{\phi_j\}_{j=1}^{rN_n} \text{ is orthonormal in } L_2(\lambda_T, [0,1])\text{]}$$

$$\leq \frac{1}{\lambda_X(K)} \left\| \sum_{j=1}^{rN_n} b_j \phi_j \right\|_\lambda^2 \quad \text{[multiplying and dividing by } \lambda_X(K)\text{]}$$

$$\leq \frac{2}{\lambda_X(K)}(\|g\|_\lambda^2 + h_0 E(a'\underline{X})^2),$$

(adding and subtracting $a'\underline{X}$ and by Assumption 2.1)

$$\leq \frac{2}{\lambda_X(K)}(1 + h_0 h_1 M^2 \lambda_{\min}^{-1})\|g\|_\lambda^2 \qquad \text{[by Assumption 2.1 and (A.10)]}.$$

From the last two displays above we conclude that (A.8) holds, with $K_1 = \frac{2}{\lambda_X(K)}(1 + h_0 h_1 M^2 \lambda_{\min}^{-1}) + h_1 \lambda_{\min}^{-1}$. This, together with (A.7), concludes the proof of this proposition. $\square$

PROOF OF THEOREM 2.1. Let $s_{I,K_n}$ be the $L_2(\mu)$ projection of $s$ onto $S_{I,K_n}$. The previous proposition and Assumptions 2.1–2.4 verify that Theorem 2.1 of Baraud (2002) can be invoked. Then

$$E\|s - \hat{s}\|_\mu^2 \leq C_2 \inf_{\mathcal{I} \times \mathcal{K}_n} \{\|s - s_{I,K_n}\|_\mu^2 + \text{pen}(I, K_n) + \vartheta_n\}.$$

Recall now that $s(\underline{x}, t) = \beta'\underline{x} + f(t)$. For any $K_n \in \mathcal{K}_n$, let $s_n(\underline{x}, t) = \beta'\underline{x} + f_{K_n}(t) \in S_{I_0,K_n}$, where we recall that $f_{K_n}$ is the $L_2(\mu_T)$ projection of $f$ onto $S_{K_n}$. Then we have

$$E\|s - \hat{s}\|_\mu^2 \leq C_2 \inf_{I_0, K_n} \{\|s - s_n\|_\mu^2 + \text{pen}(I_0, K_n) + \vartheta_n\}$$

$$\leq C_2 \inf_{K_n} \{\|f - f_{K_n}\|_{\mu_T}^2 + \text{pen}(I_0, K_n) + \vartheta_n\}. \qquad \square$$

PROOF OF COROLLARY 2.1. As an immediate consequence of Theorem 2.1 and the definition of $r_n$ (2.8), we have $\|s - \hat{s}\|_\mu^2 = O_P(r_n)$. Let now $(\underline{X}^*, T^*) \sim \mu$, with $(\underline{X}^*, T^*)$ independent of $(\underline{X}_1, T_1), \ldots, (\underline{X}_n, T_n)$. Write $E^*$ for integration with respect to $(\underline{X}^*, T^*)$ only. Notice that $E^*((\underline{X}^* - E^*(\underline{X}^*|T^*))m(T^*)) = 0$ for all bounded measurable functions $m$. Then

$$\|s - \hat{s}\|_\mu^2 = E^*(s - \hat{s})^2(\underline{X}^*, T^*)$$
$$= E^*\{(f - \hat{f})(T^*) + (\hat{\beta} - \beta)'E^*(\underline{X}^*|T^*)\}^2$$
$$+ E^*\{(\hat{\beta} - \beta)'(\underline{X}^* - E^*(\underline{X}^*|T^*))\}^2.$$

Since $\|s - \hat{s}\|_\mu^2 = O_P(r_n)$, then

(A.10) $$E^*[(\hat{\beta} - \beta)'(\underline{X}^* - E(\underline{X}^*|T^*))]^2 = O_P(r_n).$$

Since $(\underline{X}^*, T^*)$ is independent of $\hat{\beta}$, by construction we also have

(A.11) $$E^*[(\hat{\beta} - \beta)'(\underline{X}^* - E(\underline{X}^*|T^*))]^2 = (\hat{\beta} - \beta)'\Sigma(\hat{\beta} - \beta) \geq \lambda_{\min}|\hat{\beta} - \beta|_2^2.$$

Thus, from (A.10) and (A.11), and since $\lambda_{\min} > 0$, we have that for any $\beta \in \mathcal{K}$, $|\hat{\beta} - \beta|_2 = O_P(r_n^{1/2})$.

If we let $\hat{f}_1(\underline{x}) = \hat{\beta}'\underline{x}$, then we also have $\|\hat{f}_1 - f_1\|_{\mu_X} = O_P(r_n^{1/2})$. Thus

$$\|\hat{f} - f\|_{\mu_T} = \|(\hat{f} - f) + (\hat{f}_1 - f_1) - (\hat{f}_1 - f_1)\|_\mu$$
$$\leq \|\hat{s} - s\|_\mu + \|\hat{f}_1 - f_1\|_{\mu_X} = O_P(r_n^{1/2})$$

for any $f \in \mathcal{D}_{\alpha,2}(L)$.

For the empirical norm counterpart of this result, recall that we defined

$$\Gamma_n = \{\|\tilde{s}\|_\lambda \leq 2\exp(\log^2 n)\}.$$

From Baraud [(2002), proof of Theorem 1.1, page 19] we find for some constants $c_1, c_2 > 0$ that

(A.12) $\quad P(\Gamma_n^c) \leq c_1\{\exp(-2\log^2 n) + n^2 \exp(-c_2 \log^3 n)\} \to 0.$

Thus, it is enough to study the convergence rate of $\|\hat{f} - f\|_n$ on $\Gamma_n$. Notice that on $\Gamma_n$ we have $\hat{s} = \tilde{s}$, $\hat{\beta} = \tilde{\beta}$ and $\hat{f} = \tilde{f}$. Thus

$$\|\hat{f} - f\|_n^2 \mathbb{1}_{\Gamma_n} \leq 2\|\hat{s} - s\|_n^2 \mathbb{1}_{\Gamma_n} + \frac{2}{n}\sum_{i=1}^n \{(\hat{\beta} - \beta)'\underline{X}_i\}^2 \mathbb{1}_{\Gamma_n}$$

$$\leq 2\|\tilde{s} - s\|_n^2 + \frac{2}{n}\sum_{i=1}^n \{(\hat{\beta} - \beta)'\underline{X}_i\}^2 = O_P(r_n),$$

because $(\hat{\beta} - \beta)'n^{-1}\sum_{i=1}^n \underline{X}_i\underline{X}_i'(\hat{\beta} - \beta_0) = O_P(r_n)$, since $|\hat{\beta} - \beta|_2^2 = O_P(r_n)$, as above, and $n^{-1}\sum_{i=1}^n \underline{X}_i\underline{X}_i'$ converges in probability. Also $\|\tilde{s} - s\|_n^2 = O_P(r_n)$, by Corollary 3.2, page 474, in Baraud (2000), the conditions of which are verified by our assumptions and Proposition A.2. This completes the proof of this corollary. $\square$

PROOF OF COROLLARY 3.1. We evaluate now (2.8) for each penalty choice.

For the penalty term (3.4), the infimum is achieved for $K_n^* \asymp (n/\log n)^{1/2\alpha+1}$ and hence $r_n = O((\log n/n)^{2\alpha/2\alpha+1})$.

For (3.5), $r_n = O(\log n/n^{2\alpha/2\alpha+1})$ is obtained by replacing $K_n$ by $K_{n,\alpha} \asymp n^{1/2\alpha+1}$ in (2.8).

For (3.6), $r_n = O(\log n/n^{(2a+1)/(2a+2)})$ is obtained by replacing $K_n$ by $K_{n,a} \asymp n^{1/2a+2}$ in (2.8). $\square$

### A.3. The asymptotic normality of $\hat{\beta}_{q_0}$.

LEMMA A.4. *Under Assumptions 2.1–2.5, if $f \in \mathcal{D}_{\alpha,2}(L)$, we have*

(A.13) $\quad\quad\quad\quad \varepsilon'(\mathbf{Id} - \mathbf{P}_{\hat{K}_n})\mathbf{f}/\sqrt{n} \xrightarrow{P} 0,$

(A.14) $$\boldsymbol{\theta}'(\mathbf{Id} - \mathbf{P}_{\hat{K}_n})\mathbf{W}/\sqrt{n} \xrightarrow{P} 0,$$

(A.15) $$\boldsymbol{\theta}'(\mathbf{Id} - \mathbf{P}_{\hat{K}_n})\mathbf{f}/\sqrt{n} \xrightarrow{P} 0,$$

(A.16) $$\boldsymbol{\varepsilon}'\mathbf{P}_{\hat{K}_n}\mathbf{W}/\sqrt{n} \xrightarrow{P} 0 \qquad as\ n \to \infty.$$

*The above results hold uniformly over $\alpha \in (1/2, r)$ and $\gamma \in (b/4, r), b \geq 3$.*

PROOF OF THEOREM 4.1. Let $\mathbf{f} = (f(T_1), \ldots, f(T_n))'$ and $\mathbf{W} = (W_1, \ldots, W_n)'$. Recall Notation A.1, specialized now to $I = I_0$. By Lemma A.1, except for an event with probability tending to zero,

(A.17) $$\sqrt{n}(\hat{\beta}_{q_0} - \beta_{q_0}) = \sqrt{n}(\mathbf{X}'_{I_0}(\mathbf{Id} - \mathbf{P}_{\hat{K}_n})\mathbf{X}_{I_0})^{-1}\mathbf{X}'_{I_0}(\mathbf{Id} - \mathbf{P}_{\hat{K}_n})\mathbf{f} \\ + \sqrt{n}(\mathbf{X}'_{I_0}(\mathbf{Id} - \mathbf{P}_{\hat{K}_n})\mathbf{X}_{I_0})^{-1}\mathbf{X}'_{I_0}(\mathbf{Id} - \mathbf{P}_{\hat{K}_n})\mathbf{W}.$$

Since $\mathbf{X}'_{I_0} = \boldsymbol{\theta}' + \boldsymbol{\varepsilon}'$, then

(A.18) $$\sqrt{n}(\mathbf{X}'_{I_0}(\mathbf{Id} - \mathbf{P}_{\hat{K}_n})\mathbf{X}_{I_0})^{-1}\mathbf{X}'_{I_0}(\mathbf{Id} - \mathbf{P}_{\hat{K}_n})\mathbf{f} \\ = n(\mathbf{X}'_{I_0}(\mathbf{Id} - \mathbf{P}_{\hat{K}_n})\mathbf{X}_{I_0})^{-1} \\ \times (\boldsymbol{\theta}'(\mathbf{Id} - \mathbf{P}_{\hat{K}_n})\mathbf{f}/\sqrt{n} + \boldsymbol{\varepsilon}'(\mathbf{Id} - \mathbf{P}_{\hat{K}_n})\mathbf{f}/\sqrt{n}),$$

and

(A.19) $$\sqrt{n}(\mathbf{X}'_{I_0}(\mathbf{Id} - \mathbf{P}_{\hat{K}_n})\mathbf{X}_{I_0})^{-1}\mathbf{X}'_{I_0}(\mathbf{Id} - \mathbf{P}_{\hat{K}_n})\mathbf{W} \\ = n(\mathbf{X}'_{I_0}(\mathbf{Id} - \mathbf{P}_{\hat{K}_n})\mathbf{X}_{I_0})^{-1} \\ \times (\boldsymbol{\varepsilon}'\mathbf{W}/\sqrt{n} - \boldsymbol{\varepsilon}'\mathbf{P}\mathbf{W}/\sqrt{n} + \boldsymbol{\theta}'(\mathbf{Id} - \mathbf{P}_{\hat{K}_n})\mathbf{W}/\sqrt{n}).$$

By Proposition A.1 applied to $I = I_0$, $n(\mathbf{X}'_{I_0}(\mathbf{Id} - \mathbf{P}_{\hat{K}_n})\mathbf{X}_{I_0})^{-1} \xrightarrow{P} \Sigma_{q_0}^{-1}$ as $n \to \infty$. Also, notice that $\boldsymbol{\varepsilon}'\mathbf{W} = \sum_{i=1}^{n} \mathbf{D}_i$, where $\mathbf{D}_i = W_i\boldsymbol{\varepsilon}_i$ are i.i.d. vectors with $E\mathbf{D}_i = 0$ and $E(\mathbf{D}_i\mathbf{D}'_i) = \sigma^2\Sigma_{q_0}$. Then it follows from the multivariate central limit theorem that

(A.20) $$\boldsymbol{\varepsilon}'\mathbf{W}/\sqrt{n} = \sum_{i=1}^{n} \mathbf{D}_i/\sqrt{n} \xrightarrow{d} N_{q_0}(0, \sigma^2\Sigma_{q_0}).$$

Then, by (A.17)–(A.20), Lemma A.4 and Slutsky's lemma,

$$\sqrt{n}(\hat{\beta}_{q_0} - \beta_{q_0}) \xrightarrow{d} N_{q_0}(0, \sigma^2\Sigma_{q_0}^{-1}). \qquad \Box$$

PROOF OF LEMMA A.4. We begin by proving (A.13). Recall the definition of $\varepsilon_{ji}$, $j = 1, \ldots, q$, $i = 1, \ldots, n$, introduced in Appendix A.1. Then, by Assumption 2.1, for some constant $C_\varepsilon > 0$ we have $|\varepsilon_{ji}| \leq C_\varepsilon$ for all $j$ and $i$,

MODEL SELECTION IN SEMIPARAMETRIC REGRESSION 27except for a set of measure zero. Recall that $f_{K_n}$ is the $L_2(\mu_T)$ projection of $f$ onto $S_{K_n}$ and so $\mathbf{f}_{K_n} \in L_{K_n} = \{(g(T_1),\ldots,g(T_n))'|g \in S_{K_n}\}$. Then, for every $c > 0$, by Markov's inequality and (2.7),

$$P(|\boldsymbol{\varepsilon}_j'(\mathbf{Id} - \mathbf{P}_{\hat{K}_n})\mathbf{f}/\sqrt{n}| > c) \leq \sum_{k_n=A_n}^{J_n} P(|\boldsymbol{\varepsilon}_j'(\mathbf{Id} - \mathbf{P}_{2^{k_n}})\mathbf{f}/\sqrt{n}| > c)$$

$$\leq \frac{C_\varepsilon^2}{c^2} \sum_{k_n=A_n}^{J_n} E\|\mathbf{f} - \mathbf{P}_{2^{k_n}}\mathbf{f}\|_n^2 \leq \frac{C_\varepsilon^2}{c^2} \sum_{k_n=A_n}^{J_n} E\|\mathbf{f} - \mathbf{f}_{K_n}\|_n^2$$

$$= \frac{C_\varepsilon^2}{c^2} \sum_{k_n=A_n}^{J_n} \|f - f_{K_n}\|_{\mu_T}^2$$

$$\leq \frac{C_\varepsilon^2}{c^2} h_1 C(\alpha, L) \sum_{k_n=A_n}^{J_n} \frac{1}{2^{2\alpha k_n}} \to 0.$$

The second result of this lemma is very similar to (A.13). Now $\mathbf{W}$ will play the role of $\boldsymbol{\varepsilon}_j$ and $\boldsymbol{\theta}_j$ the role of $\mathbf{f}$. To emphasize more the similarity, we shall in fact prove that $\mathbf{W}'(\mathbf{Id} - \mathbf{P}_{\hat{K}_n})\boldsymbol{\theta}_j/\sqrt{n} \xrightarrow{P} 0$ for each $j \in I_0$. As above, applying Markov's inequality and using now the independence of $W$ and $(\underline{X}, T)$, we obtain

$$P(|\mathbf{W}'(\mathbf{Id} - \mathbf{P}_{\hat{K}_n})\boldsymbol{\theta}_j/\sqrt{n}| > c) \leq \sum_{k_n=A_n}^{J_n} P(|\mathbf{W}(\mathbf{Id} - \mathbf{P}_{2^{k_n}})\boldsymbol{\theta}_j/\sqrt{n}| > c)$$

$$\leq \frac{\sigma^2}{c^2} h_1 C(\gamma, A) \sum_{k_n=A_n}^{J_n} \frac{1}{2^{2\gamma k_n}} \to 0,$$

for any $c > 0$ and $\gamma > b/4$, by (2.7) applied to $\theta_j$ and recalling Assumption 2.5.

For the third result of this lemma notice that, by the same argument as in (A.6), with $\mathbf{f}$ replaced now by $\boldsymbol{\theta}_k$ we obtain

$$n^{-1/2}|\boldsymbol{\theta}_j'(\mathbf{Id} - \mathbf{P}_{\hat{K}_n})\mathbf{f}| \leq \sqrt{n}\|(\mathbf{Id} - \mathbf{P}_{\hat{K}_n})\boldsymbol{\theta}_j\|_n \|\mathbf{f} - \mathbf{P}_{\hat{K}_n}\mathbf{f}\|_n$$

$$\leq \sqrt{n}\|(\mathbf{Id} - \mathbf{P}_{\hat{K}_n})\boldsymbol{\theta}_j\|_n \|\mathbf{f} - \tilde{\mathbf{f}}\|_n,$$

where we recall now the definition of $\tilde{f}$ from Section 2.2 and that $\hat{f} = \tilde{f}\mathbb{1}_{\Gamma_n}$.

$$P(\sqrt{n}\|(\mathbf{Id} - \mathbf{P}_{\hat{K}_n})\boldsymbol{\theta}_j\|_n \|\mathbf{f} - \tilde{\mathbf{f}}\|_n \geq c)$$

$$= P(\sqrt{n}\|(\mathbf{Id} - \mathbf{P}_{\hat{K}_n})\boldsymbol{\theta}_j\|_n \|\mathbf{f} - \tilde{\mathbf{f}}\|_n \geq c, \Gamma_n)$$

$$+ P(\sqrt{n}\|(\mathbf{Id} - \mathbf{P}_{\hat{K}_n})\boldsymbol{\theta}_j\|_n \|\mathbf{f} - \tilde{\mathbf{f}}\|_n \geq c, \Gamma_n^c)$$

$$\leq P(\sqrt{n}\|(\mathbf{Id} - \mathbf{P}_{\hat{K}_n})\boldsymbol{\theta}_j\|_n \|\mathbf{f} - \hat{\mathbf{f}}\|_n \geq c) + P(\Gamma_n^c) \to 0,$$

by Corollary 3.1 and Lemma A.3, and since $P(\Gamma_n^c) \to 0$ by Theorem 2.1 in Baraud (2002).

Finally, notice that if in (A.5) we replace $\boldsymbol{\varepsilon}_k$ by $\mathbf{W}$, we obtain, up to constants, that $P(|n^{-1/2}\boldsymbol{\varepsilon}_j'\mathbf{P}_{\hat{K}_n}\mathbf{W}| > c) \leq E(\hat{K}_n)/c\sqrt{n} \to 0$, since $E(\hat{K}_n) \leq N_n \asymp (n/\log n)^{1/2}$ by (2.1). This completes the proof of this lemma. $\square$

PROOF OF COROLLARY 4.1. The proof of this corollary follows immediately by replacing $\hat{K}_n$ throughout the proofs of Theorem 4.1 and Lemma A.4 by $K_{n,\alpha}$ and $K_{n,a}$, respectively. $\square$

**A.4.**

LEMMA A.5. *Let $\mathbf{V}$ be a symmetric, semipositive definite $q \times q$ matrix. Let $\mathbf{V}^0$ be a symmetric, positive definite $q_0 \times q_0$ matrix. Let $c \in \mathbb{R}^q$ and $c_0 \in \mathbb{R}^{q_0}$, $c_0 = (c_1, \ldots, c_{q_0})$. Then the only $\mathbf{V}$ which satisfies*

$$(A.21) \qquad c_0'\mathbf{V}^0 c_0 = c'\mathbf{V} c \qquad \text{for any } c \in \mathbb{R}^q$$

*is*

$$\mathbf{V} = \begin{pmatrix} \mathbf{V}^0 & \mathbf{O}_{q_0 \times (q-q_0)} \\ \mathbf{O}_{(q-q_0) \times q_0} & \mathbf{O}_{(q-q_0) \times (q-q_0)} \end{pmatrix},$$

*where the $\mathbf{O}$ matrices have all $0$ elements.*

PROOF. Since we seek $\mathbf{V}$ symmetric and satisfying (A.21), then $\mathbf{V}$ has, in general, the structure below:

$$\mathbf{V} = \begin{pmatrix} \mathbf{V}^0 & \mathbf{A}_{q_0 \times (q-q_0)} \\ \mathbf{A}'_{(q-q_0) \times q_0} & \mathbf{B}_{(q-q_0) \times (q-q_0)} \end{pmatrix},$$

with arbitrary $\mathbf{A}$ and $\mathbf{B}$. Since we require (A.21) to hold for *any* $c \in \mathbb{R}^q$, we show that $\mathbf{A}$ and $\mathbf{B}$ are zero matrices of appropriate dimensions. First we note that for any symmetric, semipositive definite matrix $\mathbf{C}$, $x'\mathbf{C}x = 0$ for any $x$ implies $\mathbf{C} = \mathbf{O}$, with $\mathbf{O}$ being the zero matrix. Now let $c \in \mathbb{R}^q$ be arbitrary and write it as $c = (c_0, c_1)$, for $c_1 \in \mathbb{R}^{q-q_0}$. Then

$$(A.22) \qquad c'\mathbf{V}c = c_0'\mathbf{V}^0 c_0 + c_0'\mathbf{A}c_1 + c_1'\mathbf{A}'c_0 + c_1'\mathbf{B}c_1.$$

We find then $\mathbf{A}$ and $\mathbf{B}$ such that (A.22) holds, or equivalently, such that

$$(A.23) \qquad 2c_0'\mathbf{A}c_1 + c_1'\mathbf{B}c_1 = 0 \qquad \text{for any } c \in \mathbb{R}^q.$$

Since we want the above display to hold for any $c \in \mathbb{R}^q$, then it should, in particular, hold for $c = (0, c_1)$. Thus, from (A.23), $\mathbf{B}$ must satisfy

$$(A.24) \qquad c_1'\mathbf{B}c_1 = 0 \qquad \text{for any } c_1 \in \mathbb{R}^{q-q_0}.$$

MODEL SELECTION IN SEMIPARAMETRIC REGRESSION 29Now note that since $\mathbf{V}$ is semipositive definite, so is $\mathbf{B}$. This, in connection to (A.24), implies then that $\mathbf{B} = \mathbf{O}_{(q-q_0) \times (q-q_0)}$. Then we have to find $\mathbf{A}$ that satisfies $c_0' \mathbf{A} c_1 = 0$, for any $c_0 \in \mathbb{R}^{q_0}$ and $c_1 \in \mathbb{R}^{q-q_0}$. Thus, the equation must be, in particular, satisfied for the canonical basis in $\mathbb{R}^{q_0}$ and $\mathbb{R}^{q-q_0}$, respectively, which implies that $\mathbf{A} = \mathbf{O}_{q_0 \times (q-q_0)}$, which completes the proof of this lemma. $\square$

LEMMA A.6 (Rosenthal's inequality). *Let $U_1, \ldots, U_n$ be independent centered random variables with values in $\mathbb{R}$. For any $p \geq 2$, we have*

$$E \left| \sum_{i=1}^n U_i \right|^p \leq C(p) \left( E \sum_{i=1}^n |U_i|^p + \left( E \sum_{i=1}^n U_i^2 \right)^{p/2} \right),$$

*where $C(p) > 0$ depends only on $p$.*

For a proof of this inequality see, for example, Petrov (1995).

**Acknowledgments.** I am very grateful to Jon Wellner for suggesting the topic and for continuous encouragement as this work developed. I would like to thank Lucien Birgé for very helpful discussions during the Statistical Semester organized at l'Insitut Henri Poincaré, Paris, 2001, and the referees and an Associate Editor for their insightful comments.## REFERENCES

AKAIKE, H. (1974). A new look at the statistical model identification. *IEEE Trans. Automat. Control* **19** 716–723. MR423716

BARAUD, Y. (2000). Model selection for regression on a fixed design. *Probab. Theory Related Fields* **117** 467–493. MR1777129

BARAUD, Y. (2002). Model selection for regression on a random design. *ESAIM Probab. Statist.* **6** 127–146. MR1918295

BARRON, A., BIRGÉ, L. and MASSART, P. (1999). Risk bounds for model selection via penalization. *Probab. Theory Related Fields* **113** 301–413.

BICKEL, P. J., KLAASSEN, C. A. J., RITOV, Y. and WELLNER, J. A. (1993). *Efficient and Adaptive Estimation for Semiparametric Models*. John Hopkins Univ. Press. MR1679028

BIRGÉ, L. and MASSART, P. (1998). Minimum contrast estimators on sieves: Exponential bounds and rates of convergence. *Bernoulli* **4** 329–375. MR1653272

BIRGÉ, L. and MASSART, P. (2000). An adaptive compression algorithm in Besov spaces. *Constr. Approx.* **16** 1–36. MR1848840

CHEN, H. (1988). Convergence rates for parametric components in a partly linear model. *Ann. Statist.* **16** 136–146. MR924861

CHEN, H. and CHEN, K.-W. (1991). Selection of the splined variables and convergence rates in a partial spline model. *Canad. J. Statist.* **19** 323–339. MR1144149

DEVORE, R. A. and LORENTZ, G. G. (1993). *Constructive Approximation*. Springer, Berlin.

ENGLE, R., GRANGER, C., RICE, J. and WEISS, A. (1986). Semiparametric estimation of the relation between weather and electricity sales. *J. Amer. Statist. Assoc.* **81** 310–320. MR1261635

Department of Statistics
Florida State University
Tallahassee, Florida 32306-4330
USA
e-mail: bunea@stat.fsu.edu